\documentclass{book}
\usepackage[centertags]{amsmath}

\begin{document}

\title{Lectures on the mean values of functionals-An elementary introduction to infinite-dimensional probability}

\author{ Cheng-shi Liu\\Department of Mathematics\\Northeast Petroleum University\\Daqing 163318, China
\\Email:chengshiliu-68@126.com}

\maketitle

\chapter{Introduction}

The theory of functional integration namely integration on function space or infinite-dimensional integration is the important mathematical tool in probability theory and quantum and statistic physics. However, there exists a huge gap between the definitions and computations of functional integrations in general. In particular, for the norm (or more general, topology) linear space $V$, the integration on it is defined by cylinder measure depending on the linear functional space of $V$. For example, for the space $V=L^p[a,b]$, the construction of measure on it is so complicated that one almost cannot use it to compute the concrete integrations of some simple functionals  such as
\begin{equation}
f[x]=\int_a^bx(t)\mathrm{d}t,
\end{equation}
where we require $x(t)\in M=\{x(t)|0\leq x(t)\leq 1, a\leq t\leq b, x(t)\in L[a,b]\}$. It is obvious that we can expect that the average value of the function $f[x]$ on $M$ should be $\frac{b-a}{2}$. What we want is  how to obtain the result by the definition of integration on $L[a,b]$. There are several strict definitions of integration on  $L[a,b]$ through cylinder measure and its extension (see, for example, Gelfand[1], Xia[2], Gihman[3]) from which we can see that it is hopeless to use these definitions to compute the average value of $f[x]$.

Furthermore, we want to know how to solve the
following problems:

\textbf{Problem 1}. Take randomly a Lebesque integrable function $x(t)\in M$,
and let  $Y=f[x]=\int_a^bx(t)\mathrm{d}t$?. What is the probability of the event  $\frac{1}{4}<Y<\frac{1}{3}$?
What is the probability of the event  $\frac{1}{3}<Y<\frac{3}{5}$?

\textbf{Problem 2}. How to compute the mean value of the following complicated functional
\begin{equation*}
f(x)=\mathrm{e}^{\tan\int_0^1x(t)\mathrm{d}t}\sin\{\cos\int_0^1x^3(t)\mathrm{d}t\}.
\end{equation*}

Essentially, all computations of functional integrations do depend on discretization. Indeed, the cylinder measure is just a kind of discretization of measure.  Here, the difficulty is how to choice a suitable discretization in the case of Lebesque integrability.

In the lecture, we give a detailed and elementary introduction of solving these problems and other related
problems. These results show
that, in some degree, the integral in infinite dimension is more
simple than in finite dimension case. For example, the solution of
the problem 2 is
\begin{equation*}
Ef(x)=\mathrm{e}^{\tan E(\int_0^1x(t))\mathrm{d}t)}\sin\{\cos
E(\int_0^1x^3(t)\mathrm{d}t)\}=\mathrm{e}^{\tan\frac{1}{2}}\sin\{\cos
\frac{1}{4}\}.
\end{equation*}
This is just the exchange formula of the mean values of nonlinear
functionals. It is surprising since it make the computations for
some functional integrals to be very simple. For the problem 1, we
obtain the following results:
\begin{equation*}
P(\frac{1}{4}\leq\int_0^1 x(t)\mathrm{d}t\leq \frac{1}{3})=0,
\end{equation*}
\begin{equation*}
P(\frac{1}{3}\leq\int_0^1 x(t)\mathrm{d}t\leq \frac{3}{5})=1.
\end{equation*}
This is also an interesting result since it is not suitable with our
intuitions. But if we carefully study the previous first problem, these above
solutions will become the natural conclusions. In fact, we have
\begin{equation*}
E(\int_0^1 x(t)\mathrm{d}t)=\frac{1}{2},
\end{equation*}
and furthermore,
\begin{equation*}
P(\int_0^1 x(t)\mathrm{d}t=\frac{1}{2})=1.
\end{equation*}
This gives the solutions of the problem 2, and the solution of the problem 3 is also from the similar reason.

The computation of the mean values of functionals is still an important problem for mathematicians
 and physicians. In 1900's, Gateux firstly studied this problem. In 1922, Levy dealt with this topic
 in details in his book[1]. In particular, Levy obtained a famous result that now is called Levy's lemma
 which is just the concentration of measure on sphere[2-13]. Then, Wiener[14-17] introduced Wiener's
 measure to study the Brownian movement.
 In another way, Feynman's path integral became the third approach of quantum mechanics[18,19].
  The kernel of all those is the concept of functional integral and the corresponding mean values
  of functionals. Up to now, the theory of functional integration has became the mathematical foundation
  of quantum physics and then has been studied extensively[20-30]. Donsker and Varadhan studied the asymptotic evaluation of certain Markov process expectations for large time and obtained the famous lager deviation theory[31-35].

There are two basic problems in functional integral. One is the
definition of an infinite dimensional integral such as Wiener
integral and Feynman's path integral. Wiener integral has been
defined strictly[14,21,29,30], but the definition of the Feynman's path integral
 is still not satisfactory[24,25,30]. For linear topology space, one can give the definition of the infinite
 dimensional measure and hence the corresponding integral[21,22].  Another is the computation of
  functional integrals. By a direct computation, the explicit result of Gauss-type integrals
   can be obtained, and some special integrals can be reduced to Gauss type integral by transformations
   to solve[19,20,23,36-38]. For a given general functional, we can not give its explicit value. It seems that the
    computation of functional integral is more difficult than the finite dimensional case. But,
    we find that in many cases, the functional integrals is more easy to compute. Our method is
    to use the probability language and the concentration of measure. By a simple example, we give the
    basic idea and method to compute the functional integrals. Given a integral form functional
\begin{equation*}
Y=f(x)=\int_0^1 x(t)\mathrm{d}t.
\end{equation*}
Firstly, we discrete this functional as
\begin{equation*}
Y_n=\frac{1}{n}\sum_{k=1}^n x(\frac{k}{n})=\frac{1}{n}\sum_{k=1}^n x_k.
\end{equation*}
Then
\begin{equation*}
EY=\lim_{n\rightarrow \infty}EY_n=\lim_{n\rightarrow \infty}\frac{1}{n}\sum_{k=1}^n Ex_k=\frac{1}{2}.
\end{equation*}
Secondly, we compute the variance $DY$ of $Y$. By
\begin{equation*}
EY^2=E^2Y,
\end{equation*}
we give
\begin{equation*}
DY=0.
\end{equation*}
This means that the probability of $Y=EY$ is 1. It follows that a nonlinear exchange formula
\begin{equation*}
Eh(Y)=h(EY),
\end{equation*}
where $h$ can be a rather general function.

Abstractly, a functional
$f(x)$ is a function of $x$ where $x$ is an element in an
infinite-dimensional space such as $C[0,1]$. How to construct an
explicit functional by using $x(t)$ is an interesting problem. In
general, there are two basic ways to construct functional. One
method is to use the values of $x(t)$ on some points
$t_1,\cdots,t_m$ such that
\begin{equation*}
f(x)=g(x(t_1),\cdots,x(t_m)),
\end{equation*}
where $g$ is a usual function in $R^n$. Essentially, this first kind of
functional is finite dimensional functions. For example,
\begin{equation*}
f(x)=f(x_1,x_2,x_3)=x(0)+\sin(x(0.2)\exp(x(0.75)))=x_1+\sin(x_2\mathrm{e}^{x_3}),
\end{equation*}
where $x_1=x(0), x_2=x(0.2),x_3=x(0.25)$. Another method is to use
integral of $x(t)$ on some sets such that
\begin{equation*}
f(x)=\int_{I_1}\cdots\int_{I_m}g(x(t_1),\cdots,x(t_m))\mathrm{d}t_1\cdots\mathrm{d}t_m,
\end{equation*}
where $I_1,\cdots,I_m$ are subsets of the interval $[0,1]$. In
general, we take every $I_k$ be a subinterval. This second kind of
functionals is real infinite-dimensional functionals. For example,
\begin{equation*}
f(x)=\int_0^{0.3}x^2(t)\mathrm{d}t+\sin(\int_{0.1}^{0.2}x(t)\mathrm{d}t+\int_0^{0.8}\cos x(t)\mathrm{d}t).
\end{equation*}
Therefore, there are two kinds of basic elements $x(t_i)$ and
$\int_{I_i}g(x(t))\mathrm{d}t$ such that all interesting functionals
can be constructed in terms of them by addition, subtraction, multiplication, division and composition. We also call these functionals
the constructible functionals or elementary functionals.

 For the first kind of functionals, the functional integral is just the usual
integral. Thus we only consider the integral of the infinite
dimensional function. In particular, we only  consider the second
kind of functionals $f(x)$. If the domain of the functional $f$ is $M$, the
integral of $f$ on $M$ can be formally written as
\begin{equation}
\int_{M} f(x)D(x),
\end{equation}
where $D(x)$ represents formally the differential of the volume
element of $M$. But, in general, the volume $\int_{M}D(x)$ of $M$ is zero or
infinity, and the infinite-dimensional integral $\int_{M}
f(x)D(x)$ is also respectively zero or infinity. However, the mean value of functional $f$
on $M$,
\begin{equation}
Ef=\frac{\int_{M} f(x)D(x)}{\int_{M}D(x)}
\end{equation}
perhaps is finite in general. Firstly we  need a reasonable
 definition of the mean value of functional. Our approach is to use a limit procedure.
 For example, we take $M=\{x|a\leq x(t)\leq b, x(t)\in C[0,1]\}$ and $f(x)=\int_0^1g(x(t))\mathrm{d}t$,
 then we define the mean value of $f$ as
\begin{equation}
Ef=\lim_{n\rightarrow \infty}\frac{\int_a^b\cdots\int_a^b \frac{1}{n}\sum_{k=1}^ng(x_k)\mathrm{d}x_1\cdots\mathrm{d}x_n}{\int_a^b\cdots\int_a^b \mathrm{d}x_1\cdots\mathrm{d}x_n}
\end{equation}
where $x_k=x(\frac{k}{n})$. If the limitation exists and is finite or infinite, we call it the mean value of functional $f$. More general, when we take $x_k=x(t_k)$ where $t_k\in[\frac{k}{n},\frac{k+1}{n})$, if the above limitation is independent to the choose of $t_k$, we call it the mean value of $f$. In addition, we must emphasize that for difference function space, we need difference limitation procedure. In what follows, we will give concrete constructs to every case.

In the lecture, we study the functionals of integral forms
under the meaning of Riemman's and Lebesque's integrations.
Furthermore, we consider two cases including unconstraint and
constraints. We compute some kinds of functional integrals and give
the nonlinear exchange formula. Also, we discuss the problem about using
continuous functional to construct the nontrivial measurable subset
in function space such as $C[0,1]$. Our method is to consider a functional as an infinite-dimensional
random variable and then to compute the mean values and the variances of the corresponding
infinite-dimensional random variables by discretizations and limitations. In many cases, we show that the variances are zeros which mean that these functionals satisfy the property of the concentration of measure. Therefore, we give some interesting nonlinear exchange formulae on functional integrations.
In chapter 2, we study the mean
value of functional on Riemann integrable function in space $C[0,1]$.
In chapter 3, we study the mean value of functional in  function space $C[0,1]$ with the derivative constraint.
In chapter 4, we study the mean value of functional on Lebesque
integrable function space $L[0,1]$. In chapter 5, we study the mean
value of functional on the general infinite dimensional probability
space.  In chapter 6, we study the mean
value of functional on the codimension 1 subspace of function space $C[0,1]$. In chapter 7, we study the mean
value of functional on the codimension 2 subspace of  function space $C[0,1]$. In chapter 8, we give the mean values of some functionals on infinite-dimensional balls in $C[0,1]$ with 2-norm. In chapter 9, we introduce the Cauchy space and study  the
mean value of functional on it. In chapter 10, We discuss the
mean value of functional on Wiener's space.

\chapter{The mean values of functionals on $C[0,1]$}

\section{Computation of the mean value of functional on $C[0,1]$}

\subsection{Analysis method}
We consider the following functional on continuous integrable
function space $C[0,1]$,
\begin{equation*}
Y=f(x)=\int_0^1x(t)\mathrm{d}t,
\end{equation*}
where we suppose that $x(t)$ satisfies $0\leq x(t)\leq 1$, and denote $ M =\{x|0\leq x(t)\leq 1\}$.
 By discretization, we have
\begin{equation*}
Y_n=f_n(x)=\frac{1}{n}\sum_{k=0}^nx(\frac{k}{n})=\frac{1}{n}\sum_{k=0}^nx_k,
\end{equation*}
and hence
\begin{equation*}
E(Y_n)=E(f_n(x))=\frac{\int_0^1\cdots\int_0^1\frac{1}{n}\sum_{k=0}^nx_k\mathrm{d}x_1\cdots\mathrm{d}x_n}
{\int_0^1\cdots\int_0^1\mathrm{d}x_1\cdots\mathrm{d}x_n}
\end{equation*}
\begin{equation*}
=\frac{1}{n}\sum_{k=0}^n\int_0^1\cdots\int_0^1x_k\mathrm{d}x_1\cdots\mathrm{d}x_n
=\frac{1}{n}\sum_{k=0}^n\frac{1}{2}=\frac{1}{2}.
\end{equation*}
Therefore,
\begin{equation*}
EY=E(f(x))=\lim_{n\rightarrow \infty}E(Y_n)=\frac{1}{2}.
\end{equation*}

\subsection{Probability method}

We can consider $x_k$ as a random variable on [0,1]. From the probability point
of the view, we have

\begin{equation*}
E(Y_n)=\frac{1}{n}\sum_{k=0}^nE(x_k)=\frac{1}{2},
\end{equation*}
and
\begin{equation*}
E(Y)=\lim_{n\rightarrow \infty}E(Y_n)=\frac{1}{2}.
\end{equation*}

We compute the variance and density of $Y$ to give

\begin{equation*}
D(Y_n)=\frac{1}{n^2}\sum_{k=0}^nD(x_k)=\frac{1}{12n},
\end{equation*}
and then
\begin{equation*}
D(Y)=\lim_{n\rightarrow \infty}D(Y_n)=0,
\end{equation*}
\begin{equation*}
\rho_y(y)=\delta(y-\frac{1}{2}).
\end{equation*}
It means that the functional $Y=f(x)$ is almost everywhere equal to
its mean value $\frac{1}{2}$. In other words, if we choose a
function $x(t)$ randomly,  the probability of $f(x)=\frac{1}{2}$ is
1, that is
\begin{equation*}
P\{x|f(x)=\frac{1}{2}\}=1,
\end{equation*}
This is just the so-called concentration of measure.

\subsection{Characteristic function method and density function}
The characteristic function of $\frac{X_k}{n}$ is $
\frac{\mathrm{e}^{\frac{it}{n}-1}}{\frac{it}{n}}$, then the
characteristic function of $Y_n=\frac{1}{n}(X_1+\cdot+X_n)$ is
$\varphi_n(t)=(\frac{\mathrm{e}^{\frac{it}{n}-1}}{\frac{it}{n}})^n$.
Furthermore, the characteristic function of $Y$ is
\begin{equation*}
\varphi(t)=\lim_{n\rightarrow \infty}\varphi_n(t)=\lim_{n\rightarrow \infty}(\frac{\mathrm{e}^{\frac{it}{n}-1}}{\frac{it}{n}})^n=e^{\frac{\mathrm{i}t}{2}}.
\end{equation*}
It follows that the density function of $Y$ is
\begin{equation*}
\rho_Y(y)=\frac{1}{\sqrt{2\pi}}\int e^{\frac{\mathrm{i}t}{2}} e^{{-\mathrm{i}ty}}\mathrm{d}t=\delta(y-\frac{1}{2}).
\end{equation*}
This is also means that
\begin{equation*}
E(Y)=\frac{1}{2},
\end{equation*}
and
\begin{equation*}
D(Y)=0.
\end{equation*}

\section{The mean value of general functional on $C[0,1]$}

Thereafter, we denote $M=\{x|m_1\leq x(t)\leq m_2, x(t)\in R[a,b]\}$,
For the general functional
\begin{equation*}
Y=f(x)=\int_a^bg(x(t))\mathrm{d}t,
\end{equation*}
where $g$ is Riemann's integrable function, we define the mean value of $f$ on $M$ as
\begin{equation}
EY=Ef=\lim_{n\rightarrow \infty}\frac{\int_a^b\cdots\int_a^b \frac{b-a}{n}\sum_{k=1}^ng(x_k)\mathrm{d}x_1\cdots\mathrm{d}x_n}{\int_{m_1}^{m_2}\cdots\int_{m_1}^{m_2} \mathrm{d}x_1\cdots\mathrm{d}x_n},
\end{equation}
where $x_k=x(\frac{k}{n})$. By the same method, we
find the mean value of $Y$ on $M$
\begin{equation*}
EY=\frac{b-a}{m_2-m_1}\int_{m_1}^{m_2}g(x)\mathrm{d}x,
\end{equation*}
\begin{equation*}
DY=0,
\end{equation*}
\begin{equation*}
\rho_Y(y)=\delta(y-EY)
\end{equation*}

\section{The formula of the mean value of combine  functional on
$C[0,1]$}

\subsection{Probability lemmas}

Firstly, we give the following probability lemma:

\textbf{Lemma 2.1}. Assume that $A$ and $B$ are two events, and
$P(A)=P(B)=1$, then $P(AB)=1$.

\textbf{Proof}. By $P(A)=P(B)=1$, we have $P(\overline{A})=P(\overline{B})=0$. Moreover, by
$\overline{A}B\subseteq\overline{A},A\overline{B}\subseteq\overline{B},
\overline{A}\overline{B}\subseteq\overline{A}, \overline{A}\overline{B}\subseteq\overline{B}$, we have
\begin{equation*}
P(\overline{A}B)=P(A\overline{B})=P(\overline{A}\overline{B})=0.
\end{equation*}
So
\begin{equation*}
P(AB)=1-P(\overline{A}B)-P(A\overline{B})-P(\overline{A}\overline{B})=1.
\end{equation*}

According to Lemma 1 , we give the following lemmas:

\textbf{Lemma 2.2}. Assume that the densities of the random variables
$X_1,\cdots ,X_n$ are $\delta(x_1-m_1),\cdots,\delta(x_n-m_n)$. Then
the union density of $X_1,\cdots ,X_n$ is
\begin{equation*}
\rho(x_1,\cdots,x_n)=\delta(x_1-m_1)\cdots\delta(x_n-m_n).
\end{equation*}

\textbf{Lemma 2.3}. Assume that $A$ and $B$ are two events, and
$P(B)=1$, then $P(AB)=P(A)$.

\textbf{Proof}. Since $P(A\overline{B})\leq P(\overline{B})=0$, so $P(A)=P(AB)+P(A\overline{B})=P(AB)$.

According to Lemma 3, we give the following lemma:

\textbf{Lemma 2.4}. Assume that the densities of the random variables
$Y_1,\cdots ,Y_n$ are $\delta(y_1-m_1),\cdots,\delta(y_n-m_n)$, and
$\rho(x_1,\cdots,x_m)$ is the union density of $x_1,\cdots,x_m$,
then the union density of $X_1,\cdots ,X_m,Y_1,\cdots,Y_n$ is
\begin{equation*}
\rho(x_1,\cdots,x_m,y_1,\cdots,y_n)=\rho(x_1,\cdots,x_m)\delta(y_1-m_1)\cdots\delta(y_n-m_n).
\end{equation*}

\subsection{Nonlinear exchange formula}

We know that the density of the functional
 $Y=f(x)=\int_a^b g(x(t))\mathrm{d}t$ is $\rho_Y(y)=\delta(y-EY)$, so,
for the general functional $Z=h(Y)$, we have the following formula
\begin{equation*}
EZ=Eh(Y)=h(EY),
\end{equation*}
where $h(Y)$ is a general function of $Y$. In fact,
\begin{equation*}
EZ=Eh(Y)=\int h(y)\delta(y-EY)\mathrm{d}y=h(EY).
\end{equation*}
More generally, we have the following result:

\textbf{Theorem 2.1}. For $i=1,\cdots,k$, let $I_i=[a_i,b_i]$,
$Y_i=\int_{I_i} g_i(x(t))\mathrm{d}t$,
$M_i=\{x|m_{i1}\leq x(t)\leq m_{i2}, x(t)\in R[a_i,b_i]\}$,
 and  $g_i(x(t))\in R[a_i,b)i]$ with respect $t$, $g_i(x)\in R[m_{i1},m_{i2}]$ with respect to $x$.  Furthermore, let $h(y_1,\cdots,y_k)$
be a general function. Then we have
\begin{equation*}
Eh(Y_1,\cdots,Y_k)=h(EY_1,\cdots,EY_k),
\end{equation*}
where
\begin{equation*}
EY_i=\frac{b_i-a_i}{m_{i2}-m_{i1}}\int_{m_{i1}}^{m_{i2}}g_i(x)\mathrm{d}x.
\end{equation*}

\textbf{Proof.} Since every $EY_i$ is the mean value of $Y_i$ on $M_i$ which definition is given by
\begin{equation}
EY_i=\lim_{n\rightarrow \infty}\frac{\int_{a_i}^{b_i}\cdots\int_{a_i}^{b_i} \frac{b_i-a_i}{n}\sum_{k=1}^ng(x_k)\mathrm{d}x_1\cdots\mathrm{d}x_n}{\int_{m_{i1}}^{m_{i2}}\cdots\int_{m_{i1}}^{m_{i2}} \mathrm{d}x_1\cdots\mathrm{d}x_n},
\end{equation}
where $x_k=x(\frac{k}{n})$, we similarly have
\begin{equation}
DY_i=0,
\end{equation}
and the density of $Y_i$ is $\delta(y_i-EY_i)$, and then
\begin{equation*}
Eh(Y_1,\cdots,Y_k)=\int\cdots\int h(y_1,\cdots,y_k)\delta(y_1-EY_1)\cdots \delta(y_k-EY_k)\mathrm{d}y_1\cdots\mathrm{d}y_k
\end{equation*}
\begin{equation*}
=h(EY_1,\cdots,EY_k).
\end{equation*}

\textbf{Example 2.1}. For $Y=f(x)=\sin(\int_0^{\frac{1}{2}}x^3(t)\mathrm{d}t)\mathrm{e}^{(\int_0^1x(t)\mathrm{d}t)^2}$, $M=\{x|0\leq x(t)\leq 1\}$, then its mean value on $M$ is
\begin{equation*}
EY=\mathrm{e}^{\frac{1}{4}}\sin\frac{1}{8}.
\end{equation*}

By lemma 2.4, we have the following theorem.

\textbf{Theorem 2.2}. For $i=1,\cdots,k$, let $I_i=[a_i,b_i]$,
$Y_i=\int_{I_i} g_i(x(t))\mathrm{d}t$,
$M_i=\{x|m_{i1}\leq x(t)\leq m_{i2}, x(t)\in R[a_i,b_i]\}$,
 and  $g_i(x(t))\in R[a_i,b)i]$ with respect $t$, $g_i(x)\in R[m_{i1},m_{i2}]$ with respect to $x$.  Furthermore, let $h(x_1,\cdots,x_h; y_1,\cdots,y_k)$
be a general function. Then we have
\begin{equation*}
Eh(X_1,\cdots,X_m,Y_1,\cdots,Y_k)=Eh(X_1,\cdots,X_m,EY_1,\cdots,EY_k).
\end{equation*}

\chapter{The averages of functionals on $L[0,1]$}

 In this section, we study the mean values of
functionals  on the Lebesque's integrable function space
$L[0,1]$.

\section{The mean value of functional $Y=f(x)=\int_0^1x(t)\mathrm{d}t$}
We consider the following functional on the subset $M=\{x|0\leq x(t)\leq1, x(t)\in L[0,1]\}$ on Lebesque's integrable function
space $L[0,1]$,
\begin{equation*}
Y=f(x)=\int_0^1x(t)\mathrm{d}t.
\end{equation*}
By discretization, we have
\begin{equation*}
Y_n=\sum_{k=1}^n\frac{k}{n}m(E_k),
\end{equation*}
where $E_k=\{t|\frac{k-1}{n}<x(t)\leq\frac{k}{n}\}$, and $m(E_k)$ is the Lebesque measure of $E_k$ satisfying
\begin{equation*}
\sum_{k=1}^n m(E_k)=1.
\end{equation*}
Denote $z_k=m(E_k)$, and consider $z_k$ as random variable. Then we have
\begin{equation*}
Y_n=\sum_{k=1}^n\frac{k}{n}z_k,
\end{equation*}
with
\begin{equation*}
\sum_{k=1}^n z_k=1,
\end{equation*}
and $0\leq z_k\leq1$. We define the mean value of $Y$ on $M$ as follows
\begin{equation}
EY=\lim_{n\rightarrow \infty}\frac{\int_{M_0} \sum_{k=1}^n\frac{k}{n}z_k\mathrm{d}V}{\int_{M_0}\mathrm{d}V}
\end{equation}
where $M_0=\{(z_1,\cdots,z_n)|z_1+\cdots+z_n=1, z_1\geq0,\cdots,z_n\geq0\}$, and $\mathrm{d}V$ is the volume element of $M_0$.

Therefore, because the mean value of $Y_n$ is
\begin{equation*}
EY_n=\sum_{k=1}^n\frac{k}{n}E(z_k)=E(z_1)\sum_{k=1}^n\frac{k}{n}=\frac{n+1}{2}E(z_1),
\end{equation*}
we only need to compute $E(z_1)$. There are two methods to get its value $E(z_1)=\frac{1}{n}$. One simple method is from $E(z_1)=\cdots=E(z_n)$  and $\sum_{k=1}^n z_k=1$. Another method is to direct compute the corresponding integrals.
Denote  $M_1=\{(z_1,\cdots,z_n)|\sum_{k=1}^{n-1} z_k\leq1,
 z_k\geq 0\}$. We need the following results.

\textbf{Lemma 3.1}.
\begin{equation*}
\int_{M_0}\mathrm{d}V=\frac{\sqrt n}{(n-1)!},
\end{equation*}
\begin{equation*}
\int_{M_0}z_1\mathrm{d}V=\frac{\sqrt n}{n!},
\end{equation*}
\begin{equation*}
\int_{M_0}z_1z_2\mathrm{d}V=\frac{\sqrt n}{(n+1)!}.
\end{equation*}
\begin{equation*}
\int_{M_0}z^2_1\mathrm{d}V=\frac{2\sqrt n}{(n+1)!}.
\end{equation*}
\textbf{Proof.} By direct computation, we have
\begin{equation*}
\int_{M_0}\mathrm{d}V=\int_{M_1}\sqrt{1+(\frac{\partial z_n}{\partial z_1})^2+\cdots+(\frac{\partial z_n}{\partial z_{n-1}})^2}\mathrm{d}z_1\cdots\mathrm{d}z_{n-1}
\end{equation*}
\begin{equation*}
=\int_{M_1}\sqrt{n}\mathrm{d}z_1\cdots\mathrm{d}z_{n-1}=\frac{\sqrt n}{(n-1)!};
\end{equation*}

\begin{equation*}
\int_{M_0} z_1\mathrm{d}V=\int_{M_1}z_1\sqrt{1+(\frac{\partial z_n}{\partial z_1})^2+\cdots+(\frac{\partial z_n}{\partial z_{n-1}})^2}\mathrm{d}z_1\cdots\mathrm{d}z_{n-1}
\end{equation*}
\begin{equation*}
=\int_{M_1}\sqrt{n}z_1\mathrm{d}z_1\cdots\mathrm{d}z_{n-1}=\frac{\sqrt n}{(n-2)!}\int_0^1z_1(1-z_1)^{n-2}\mathrm{d}z_1=\frac{\sqrt n}{n!}.
\end{equation*}
\begin{equation*}
\int_{M_0}z_1z_2\mathrm{d}V=\int_{M_1}z_1z_2\sqrt{1+(\frac{\partial z_n}{\partial z_1})^2+\cdots+(\frac{\partial z_n}{\partial z_{n-1}})^2}\mathrm{d}z_1\cdots\mathrm{d}z_{n-1}
\end{equation*}
\begin{equation*}
=\frac{\sqrt n}{(n-3)!}\int_0^1z_1\mathrm{d}z_1\int_0^{1-z_1}\mathrm{d}z_2(1-z_1-z_2)^{n-3}=\frac{\sqrt n}{(n+1)!};
\end{equation*}

\begin{equation*}
\int_{M_0} z^2_1\mathrm{d}V=\int_{M_1}z^2_1\sqrt{1+(\frac{\partial z_n}{\partial z_1})^2+\cdots+(\frac{\partial z_n}{\partial z_{n-1}})^2}\mathrm{d}z_1\cdots\mathrm{d}z_{n-1}
\end{equation*}
\begin{equation*}
=\int_{M_1}\sqrt{n}z^2_1\mathrm{d}z_1\cdots\mathrm{d}z_{n-1}=\frac{\sqrt n}{(n-2)!}\int_0^1z^2_1(1-z_1)^{n-2}\mathrm{d}z_1=\frac{2\sqrt n}{(n+1)!}.
\end{equation*}

By the above lemmas, we have

\textbf{Lemma 3.2.}
\begin{equation*}
E(z_1)=\frac{\sqrt n}{n!}/\frac{\sqrt n}{(n-1)!}=\frac{1}{n},
\end{equation*}
\begin{equation*}
E(z_1z_2)=\frac{\sqrt n}{(n+1)!}/\frac{\sqrt n}{(n-1)!}=\frac{1}{n(n+1)},
\end{equation*}
\begin{equation*}
E(z^2_1)=\frac{2\sqrt n}{(n+1)!}/\frac{\sqrt n}{(n-1)!}=\frac{2}{n(n+1)}.
\end{equation*}

By these results, we can give $E(Y)$ and $D(Y)$.

\textbf{Theorem 3.1.} For the above functional $Y$ and $M$, we have
\begin{equation*}
E(Y)=\lim_{n\rightarrow +\infty}E(Y_n)=\frac{1}{2},
\end{equation*}
\begin{equation*}
D(Y)=0.
\end{equation*}
Further, for a general function $h(Y)$, we have
\begin{equation*}
E(h(Y))=h(\frac{1}{2}).
\end{equation*}

\textbf{Proof}. In fact, we have
\begin{equation*}
E(Y_n)=\frac{n+1}{2n},
\end{equation*}
and
\begin{equation*}
E(Y)=\lim_{n\rightarrow +\infty}E(Y_n)=\frac{1}{2}.
\end{equation*}
Now we compute the variances of $Y_n$ and $Y$. Firstly, we have
\begin{equation*}
EY^2_n=\sum_{j=1}^n\sum_{k=1}^n\frac{jk}{n^2}E(z_kz_j)=\frac{1}{n^2}\sum_{k=1}^n k^2E(z_k^2)+\frac{1}{n^2}\sum_{k\neq j}kjE(z_kz_j)
\end{equation*}
\begin{equation*}
=\frac{1}{n^3(n+1)}(\sum_{k=1}^n 2k^2+\sum_{j\neq k}jk)=\frac{1}{n^3(n+1)}(\sum_{k=1}^n k^2+\sum_{j, k}jk)
\end{equation*}
\begin{equation*}
=\frac{1}{n^3(n+1)}\{\frac{n^2(n+1)^2}{4}+\frac{n(n+1)(2n+1)}{6}\},
\end{equation*}
and
\begin{equation*}
E(Y^2)=\lim_{n\rightarrow +\infty}\frac{1}{n^3(n+1)}\{\frac{n^2(n+1)^2}{4}+\frac{n(n+1)(2n+1)}{6}\}=\frac{1}{4}.
\end{equation*}
So
\begin{equation*}
D(Y)=0.
\end{equation*}
This means that the density of $Y$ is $\rho(y)=\delta(y-\frac{1}{2})$. Therefore, we have
\begin{equation*}
E(h(Y))=h(\frac{1}{2}),
\end{equation*}
where $h(Y)$ is a general function. The proof is completed.

For example, $E(\sin(\int_0^1x(t)\mathrm{d}t)=\sin\frac{1}{2}$.

\section{The mean value of functional $Y=\int_0^1g(x(t))\mathrm{d}t$ on
$L[0,1]$}

\textbf{Theorem 3.2}. Assume that $g(x)$ is a piecewise monotonic  continuous function of $x$ on $[0,1]$  and
$x(t)$ is a measurable function  on $[0,1]$,  $M=\{x|a\leq x(t)\leq b, x(t)\in L[0,1]\}$. For the functional
\begin{equation*}
Y=\int_0^1g(x(t))\mathrm{d}t,
\end{equation*}
we have
\begin{equation*}
E(Y)=\lim_{n\rightarrow +\infty}E(Y_n)=\int_0^1g(x)\mathrm{d}x,
\end{equation*}
and
\begin{equation*}
D(Y)=0.
\end{equation*}
Further, we have
\begin{equation*}
E(h(Y))=h(EY),
\end{equation*}
where $h(Y)$ is a general function.

\textbf{Proof}. We first assume $g(x)$ is a monotonic function on whole interval $[0,1]$. Then the discretization gives
\begin{equation}
Y_n=\sum_{k=1}^ng(\frac{k}{n})m(E_k),
\end{equation}
where for $g(x)$ being an increasing function we have $E_k=\{t|g(\frac{k-1}{n})<x(t)\leq g(\frac{k}{n})\}=\{t|\frac{k-1}{n}<x(t)\leq \frac{k}{n}\}$, for $g(x)$ being a  decreasing function we have $E_k=\{t|g(\frac{k}{n})<x(t)\leq g(\frac{k-1}{n})\}=\{t|\frac{k-1}{n}<x(t)\leq \frac{k}{n}\}$, and $m(E_k)$ is the Lebesque measure of $E_k$ satisfying
\begin{equation*}
\sum_{k=1}^n m(E_k)=1.
\end{equation*}
Denote $z_k=m(E_k)$, and consider $z_k$ as the random variable. Then
we have
\begin{equation}
Y_n=\sum_{k=1}^ng(\frac{k}{n})z_k,
\end{equation}
with
\begin{equation*}
\sum_{k=1}^n z_k=1,
\end{equation*}
and $0\leq z_k\leq1$. Therefore, the mean value of $Y_n$ is
\begin{equation*}
EY_n=\sum_{k=1}^ng(\frac{k}{n})E(z_k)=\frac{1}{n}\sum_{k=1}^ng(\frac{k}{n}).
\end{equation*}
So, we have
\begin{equation*}
E(Y)=\lim_{n\rightarrow +\infty}E(Y_n)=\int_0^1g(x)\mathrm{d}x.
\end{equation*}
Furthermore, we give
\begin{equation*}
EY^2_n=\sum_{k=1}^n\sum_{j=1}^n g(\frac{k}{n})g(\frac{j}{n})E(z_kz_j)
\end{equation*}
\begin{equation*}
=\sum_{k=1}^n g^2(\frac{k}{n})E(z^2_k)+\sum_{k\neq j} g(\frac{k}{n})g(\frac{j}{n})E(z_kz_j)
\end{equation*}
\begin{equation*}
=\frac{2}{n(n+1)}\sum_{k=1}^n g^2(\frac{k}{n})+\frac{1}{n(n+1)}\sum_{k\neq j} g(\frac{k}{n})g(\frac{j}{n}),
\end{equation*}
and then
\begin{equation*}
E(Y^2)=\lim_{n\rightarrow +\infty}E(Y^2_n)=(\int_0^1g(x)\mathrm{d}x)^2.
\end{equation*}
So
\begin{equation*}
D(Y)=0.
\end{equation*}
This means that the density of $Y$ is $\rho(y)=\delta(y-EY)$. Therefore, we have
\begin{equation*}
E(h(Y))=h(EY),
\end{equation*}
where $h(Y)$ is a general function. For the case of $g$ being a piecewise monotonic function, we can consider respectively each subinterval on which $g$ is monotonic, and then combine these results to give the conclusion.  The proof is completed.

For example, we take $M=\{x|0\leq x(t)\leq 1, x(t)\in L[0,1]\}$ and have $E(\sin(\int_0^1x^2(t)\mathrm{d}t)=\sin\frac{1}{3}$.

\textbf{Remark 3.1}. Since $g(x)$ is a continuous function, for $n$ being enough large, $|g(\frac{k-1}{n})-g(\frac{k}{n})|$ will be enough small, so the discretization (9)  is reasonable from the definition of Lebesque integral.

\chapter{The mean values of functionals on $C[0,1]$ with derivative
constraints}

Consider the arc length of a curve $x(t)$
defined on [0,1],
\begin{equation*}
l(x)=\int_0^1 \sqrt{1+(x'(t))^2}\mathrm{d}t.
\end{equation*}
We want to compute the mean value of the arc length on some set $M$ of curves. For example, we take
 $M=\{0\leq x(t)\leq1, x(t)\in C^1[0,1]\}$. Although $x(t)$ is bounded, $x'(t)$ perhaps is unbounded. So the mean value of
 arc length will be infinity. We consider a more simple example to give an infinity result.
 Take a functional
 \begin{equation*}
 Y=\int_0^1(x'(t))^2\mathrm{d}t,
 \end{equation*}
where $x'(t)$ is continuous function. Correspondingly, we take the disretization of $Y$ as
\begin{equation*}
Y_n=\frac{1}{n}\sum_{k=1}^n (\frac{x(\frac{k}{n})-x(\frac{k-1}{n})}{1/n})^2=n\sum_{k=1}^n(x^2_k+x^2_{k-1}-2x_kx_{k-1}).
\end{equation*}
So,
\begin{equation*}
EY_n=n\sum_{k=1}^n(Ex^2_k+Ex^2_{k-1}-2Ex_kEx_{k-1})=\frac{n^2}{6},
\end{equation*}
and then
\begin{equation*}
EY=\lim_{n\rightarrow +\infty}EY_n=\lim_{n\rightarrow +\infty}\frac{n^2}{6}=+\infty .
\end{equation*}
In order to obtain a finite value $EY<\infty$, we need to add some
constraint on $x'(t)$. For example, take $M=\{0\leq x'(t)\leq 1\}$.
We compute the mean value of $Y$ in $M$ to give
\begin{equation*}
EY=E(\int_0^1 z(t)\mathrm{d}t)=\frac{1}{3},
\end{equation*}
where $z(t)=x'(t)$.

However, under the constraint on derivative function $x'(t)$, to
compute the mean value of a functional constructed by $x(t)$, we
need a little trick. In next section, we will give the details.

\section{The mean value of the functional $Y=\int_0^1 x(t)\mathrm{d}t$
with derivative constraint}

\textbf{Theorem 4.1}. Take the constraint $M=\{0\leq x'(t)\leq1,x(0)=0\}$. For  the
mean value of the functional $Y=\int_0^1 x(t)\mathrm{d}t$ on $M$, we have
\begin{equation*}
EY=\frac{1}{4},
\end{equation*}
and
\begin{equation*}
D(Y)=0.
\end{equation*}

\textbf{Proof}. We first transform $Y$ to a functional of $x'(t)$. In
fact, by Newton-Leibnitz's formula, we have
\begin{equation*}
x(t)=\int_0^t x'(s)\mathrm{d}s+x(0)=\int_0^t x'(s)\mathrm{d}s,
\end{equation*}
on $M$. Then $Y$ can be rewritten as
\begin{equation*}
Y=\int_0^1\int_0^t x'(s)\mathrm{d}s\mathrm{d}t=\int_0^1x'(t)\mathrm{d}t-\int_0^1tx'(t)\mathrm{d}t.
\end{equation*}
Let
$z(t)=x'(t)$,$Y_1=\int_0^1x'(t)\mathrm{d}t=\int_0^1z(t)\mathrm{d}t$
and $Z=\int_0^1tx'(t)\mathrm{d}t=\int_0^1tz(t)\mathrm{d}t$. Then we
get $EY_1=\frac{1}{2}$ on $M$. Indeed, we only need to compute $EZ$.
We have
\begin{equation*}
Z_n=\frac{1}{n}\sum_{k=1}^n\frac{k}{n}z_k=\frac{1}{n^2}\sum_{k=1}^nkz_k,
\end{equation*}
and then
\begin{equation*}
EZ_n=\frac{1}{n^2}\sum_{k=1}^nkEz_k=\frac{n+1}{4n}.
\end{equation*}
Therefore,
\begin{equation*}
EZ=\lim_{n\rightarrow +\infty}EZ_n=\lim_{n\rightarrow +\infty}\frac{n+1}{4n}=\frac{1}{4}.
\end{equation*}
In order to get the variance $D(Z)$ of $Z$, we first compute $EZ^2$.
We have
\begin{equation*}
Z^2_n=\frac{1}{n^4}\sum_{j=1}^n\sum_{k=1}^njkz_jz_k,
\end{equation*}
\begin{equation*}
EZ^2_n=\frac{1}{n^4}\sum_{j=1}^n\sum_{k=1}^njkE(z_jz_k)
\end{equation*}
\begin{equation*}
=\frac{1}{n^4}\sum_{k=1}^nk^2E(z^2_k)+\frac{1}{n^4}\sum_{j\neq k}jkEz_jEz_k
\end{equation*}
\begin{equation*}
=\frac{1}{3n^4}\sum_{k=1}^nk^2+\frac{1}{4n^4}\sum_{j\neq k}jk
\end{equation*}
\begin{equation*}
=\frac{(n+1)(2n+1)}{4n^3}+\frac{(n+1)^2}{16n^2}.
\end{equation*}
So
\begin{equation*}
EZ^2=\lim_{n\rightarrow +\infty}\{\frac{(n+1)(2n+1)}{4n^3}+\frac{(n+1)^2}{16n^2}\}=\frac{1}{16} .
\end{equation*}
Hence, we have
\begin{equation*}
D(Z)=0.
\end{equation*}
This means that the density of $Z$ is $\delta(z-\frac{1}{4})$. We know also the density of $Y_1$ is $\delta(y-\frac{1}{2})$. It follows that
\begin{equation*}
EY=EY_1-EZ=\frac{1}{2}-\frac{1}{4}=\frac{1}{4},
\end{equation*}
and
\begin{equation*}
D(Y)=0.
\end{equation*}
The proof is completed.

\section{The mean value of the functional $Y=\int_0^1 x^2(t)\mathrm{d}t$
with derivative constraint}

\textbf{Theorem 4.2}. Take the constraint $M=\{0\leq
x'(t)\leq1,x(0)=0\}$. For the mean value of the functional $Y=\int_0^1 x^2(t)\mathrm{d}t$, we have
\begin{equation*}
EY=EY_1-EZ=\frac{1}{4}-2\times\frac{1}{12}=\frac{1}{12},
\end{equation*}
and
\begin{equation*}
D(Y)=0.
\end{equation*}

\textbf{Proof}. $Y$ can be rewritten as
\begin{equation*}
Y=\int_0^1(\int_0^t x'(s)\mathrm{d}s)^2\mathrm{d}t=(\int_0^1x'(t)\mathrm{d}t)^2-2\int_0^1tx'(t)\int_0^tx'(s)\mathrm{d}s\mathrm{d}t.
\end{equation*}
Let $z(t)=x'(t)$,$Y_1=(\int_0^1x'(t)\mathrm{d}t)^2=(\int_0^1z(t)\mathrm{d}t)^2$ and $Z=\int_0^1tx'(t)\int_0^tx'(s)\mathrm{d}s\mathrm{d}t=\int_0^1tz(t)\int_0^tz(s)\mathrm{d}s\mathrm{d}t$. Then $EY_1=\frac{1}{4}$ on $M$. We only need to compute $EZ$.
We have
\begin{equation*}
Z_n=\frac{1}{n}\sum_{k=1}^n\frac{k}{n}z_k\frac{1}{n}\sum_{j=0}^{k-1}z_j=\frac{1}{n^3}\sum_{k=1}^nk\sum_{j=0}^{k-1}z_jz_k,
\end{equation*}
and then
\begin{equation*}
EZ_n=\frac{1}{n^3}\sum_{k=1}^nk\sum_{j=0}^{k-1}Ez_jEz_k=\frac{(n+1)(2n+1)}{24n^2}.
\end{equation*}
Therefore,
\begin{equation*}
EZ=\lim_{n\rightarrow +\infty}EZ_n=\lim_{n\rightarrow +\infty}\frac{(n+1)(2n+1)}{24n^2}=\frac{1}{12}.
\end{equation*}
It follows that
\begin{equation*}
EY=EY_1-2EZ=\frac{1}{4}-2\times\frac{1}{12}=\frac{1}{12}.
\end{equation*}
In order to get the variance $D(Z)$ of $Z$, we first compute $EZ^2$.
We have
\begin{equation*}
Z^2_n=\frac{1}{n^6}\sum_{k=1}^nk\sum_{j=0}^{k-1}z_jz_k\sum_{h=1}^nh\sum_{i=0}^{h-1}z_iz_h
\end{equation*}
\begin{equation*}
=\frac{1}{n^6}\sum_{k=1}^n \sum_{h=1}^nkhz_kz_h\sum_{j=0}^{k-1}\sum_{i=0}^{h-1}z_i z_j
\end{equation*}
\begin{equation*}
=\frac{1}{n^6}\sum_{k=1}^nk^2z^2_k\sum_{j=0}^{k-1}\sum_{i=0}^{k-1}z_i z_j+\frac{2}{n^6}\sum_{h=2}^n \sum_{k=1}^{h-1}khz_kz_h(\sum_{j=0}^{k-1}\sum_{i=0}^{k-1}z_i z_j+\sum_{j=0}^{k-1}\sum_{i=k}^{h-1}z_i z_j)
\end{equation*}
\begin{equation*}
EZ^2_n=\frac{1}{n^6}\sum_{k=1}^nk^2Ez^2_k\sum_{j=0}^{k-1}\sum_{i=0}^{k-1}E(z_iz_j)+\
\end{equation*}
\begin{equation*}
+\frac{2}{n^6}\sum_{h=2}^n \sum_{k=1}^{h-1}khEz_kEz_h(\sum_{j=0}^{k-1}\sum_{i=0}^{k-1}E(z_iz_j)+\sum_{j=0}^{k-1}\sum_{i=k}^{h-1}Ez_iEz_j).
\end{equation*}
By $\sum_{j=0}^{k-1}\sum_{i=0}^{k-1}E(z_iz_j)=\frac{k^2}{4}+\frac{k}{12}$, we have
\begin{equation*}
Ez^2_n=\frac{1}{3n^6}\sum_{k=1}^nk^2(\frac{k^2}{4}+\frac{k}{12})+\frac{2}{n^6}\sum_{h=2}^n
\sum_{k=1}^{h-1}\frac{1}{4}kh(\frac{k^2}{4}+\frac{k}{12}+\frac{1}{4}(h-k-1)k)
\end{equation*}
\begin{equation*}
=\frac{1}{48n^6}\sum_{h=2}^nh^5+low order terms
\end{equation*}
So
\begin{equation*}
EZ^2=\lim_{n\rightarrow +\infty}\{\frac{1}{48n^6}\frac{n^2(n+1)^2(2n^2+2n-1)}{12}+low order terms\}=\frac{1}{12^2} .
\end{equation*}
Hence, we have
\begin{equation*}
D(Z)=0.
\end{equation*}
This means that the density of $Z$ is $\delta(z-\frac{1}{12})$. We know also the density of $Y_1$ is $\delta(y-\frac{1}{4})$. It follows that
\begin{equation*}
EY=EY_1-EZ=\frac{1}{4}-2\times\frac{1}{12}=\frac{1}{12},
\end{equation*}
and
\begin{equation*}
D(Y)=0.
\end{equation*}
The proof is completed.

\section{The mean value of the functional $Y=\int_0^1
g(x(t))\mathrm{d}t$ with derivative constraint}

\textbf{Theorem 4.3}. Take the constraint
$M=\{0\leq x'(t)\leq1,x(0)=0\}$, and let $g(x)$ is differentiable. For the mean value of the functional $Y=\int_0^1
g(x(t))\mathrm{d}t$, we have
\begin{equation*}
EY=g(\frac{1}{2})-\frac{1}{2}\int_0^1tg'(\frac{t}{2})\mathrm{d}t
=\frac{1}{4}-\frac{1}{2}\int_0^1t^2\mathrm{d}t=\frac{1}{4}-\frac{1}{6}=\frac{1}{12},
\end{equation*}
\begin{equation*}
D(Y)=0.
\end{equation*}

\textbf{Proof}. Firstly, $Y$ can be rewritten as
\begin{equation*}
Y=\int_0^1g(\int_0^t x'(s)\mathrm{d}s))\mathrm{d}t=g(\int_0^1x'(t)\mathrm{d}t)-\int_0^1tx'(t)g'(\int_0^tx'(s)\mathrm{d}s)\mathrm{d}t.
\end{equation*}
Let $z(t)=x'(t)$,$Y_1=g(\int_0^1x'(t)\mathrm{d}t)=g(\int_0^1z(t)\mathrm{d}t)$ and $Z=\int_0^1tx'(t)g'(\int_0^tx'(s)\mathrm{d}s)\mathrm{d}t=\int_0^1tz(t)g'(\int_0^tz(s)\mathrm{d}s)\mathrm{d}t$. Then $EY_1=g(\frac{1}{2})$ on $M$. We only need to compute $EZ$. For the purpose we use the previous formula to give the corresponding results. Indeed, according to $z(t)$  being independent with $z(s)$ for $t\neq s$,  we have
\begin{equation*}
EZ=\int_0^1tE(z(t))g'(\int_0^tE(z(s))\mathrm{d}s)\mathrm{d}t=\frac{1}{2}\int_0^1tg'(\frac{t}{2})\mathrm{d}t.
\end{equation*}
It follows that
\begin{equation*}
EY=EY_1-EZ=g(\frac{1}{2})-\frac{1}{2}\int_0^1tg'(\frac{t}{2})\mathrm{d}t.
\end{equation*}
For example,  taking $g(x)=x^2$, then
\begin{equation*}
EY=g(\frac{1}{2})-\frac{1}{2}\int_0^1tg'(\frac{t}{2})\mathrm{d}t
=\frac{1}{4}-\frac{1}{2}\int_0^1t^2\mathrm{d}t=\frac{1}{4}-\frac{1}{6}=\frac{1}{12}.
\end{equation*}
By the same reason, we know that $DY=0$. The proof is completed.

\textbf{Remark 4.1}. The theorems 4.1 and 4.2 are the special cases oh the theorem 4.3.

\section{The mean value of the arc length on $M=\{x(t)|0\leq x'(t)\leq
1\}$}

\textbf{Theorem 4.4}. Take $M=\{x(t)|0\leq x'(t)\leq
1\}$.  For the mean value of the arc lengths $Y$ of $x(t)$ on $M$, we have
\begin{equation*}
EY=\frac{\sqrt2}{2}+\frac{1}{2}\ln(1+\sqrt2),
\end{equation*}
\begin{equation*}
D(Y)=0.
\end{equation*}

\textbf{Proof}. The arc length formula is
\begin{equation*}
Y=l(x)=\int_0^1\sqrt{1+(x'(t))^2}\mathrm{d}t.
\end{equation*}
Taking $z(t)=x'(t)$, the above formula becomes
\begin{equation*}
Y=l(z)=\int_0^1\sqrt{1+z^2(t)}\mathrm{d}t.
\end{equation*}
Correspondingly, we have
\begin{equation*}
Y_n=\frac{1}{n}\sum_{k=1}^n\sqrt{1+z^2_k}.
\end{equation*}
So
\begin{equation*}
EY_n=\frac{1}{n}\sum_{k=1}^nE\sqrt{1+z^2_k}=E\sqrt{1+z^2_1}.
\end{equation*}
Moreover,
\begin{equation*}
E\sqrt{1+z^2_1}=\int_0^1\sqrt{1+z^2_1}\mathrm{d}z_1=\frac{\sqrt2}{2}+\frac{1}{2}\ln(1+\sqrt2).
\end{equation*}
Hence, the mean value of arc length on M is
\begin{equation*}
EY=\int_0^1\sqrt{1+z^2_1}\mathrm{d}z_1=\frac{\sqrt2}{2}+\frac{1}{2}\ln(1+\sqrt2),
\end{equation*}
and $DY=0$. This also means that the density of $Y$ is $\rho(y)=\delta(y-\frac{\sqrt2}{2}-\frac{1}{2}\ln(1+\sqrt2))$. The proof is completed.

\chapter{The mean values of functionals on infinite-dimensional space with
general probability measure}

In previous sections, we still suppose that the random variable is uniform distribution on $[0,1]$. Now we consider the general random variable.  For the following functional on $R[0,1]$,
\begin{equation*}
Y=f(x)=\int_0^1x(t)\mathrm{d}t,
\end{equation*}
we suppose that random variable $x(t)$ satisfies $-\infty\leq x(t)\leq +\infty$ and its density is $\rho(x)$. In sections 8 and 9, we will see the concrete examples.
By
\begin{equation*}
Y_n=f_n(x)=\frac{1}{n}\sum_{k=0}^nx(\frac{k}{n})=\frac{1}{n}\sum_{k=0}^nx_k,
\end{equation*}
we have
\begin{equation*}
EY_n=\frac{1}{n}\sum_{k=0}^nEx_k=Ex_1=\int x\rho(x)\mathrm{d}x,
\end{equation*}
and hence
\begin{equation*}
EY=E(f(x))=\lim_{n\rightarrow \infty}E(Y_n)=\int x\rho(x)\mathrm{d}x.
\end{equation*}
By
\begin{equation*}
D(Y_n)=\frac{1}{n^2}\sum_{k=0}^nD(x_k)=\frac{1}{n}D(x_1),
\end{equation*}
if $Dx_1<+\infty$, we have
\begin{equation*}
D(Y)=\lim_{n\rightarrow \infty}D(Y_n)=0,
\end{equation*}
and then
\begin{equation*}
\rho_y(y)=\delta(y-EY),
\end{equation*}
It means that the functional $Y=f(x)$ is almost everywhere equal to
its mean value $EY$. In other words, if we choose a function $x(t)$
randomly,  the probability of $f(x)=EY$ is 1, that is
\begin{equation*}
P(x|f(x)=EY)=1,
\end{equation*}
This is just the concentration of measure.

We consider the following general function
\begin{equation*}
Y=f(x)=\int_a^b g(x(t))\mathrm{d}t,
\end{equation*}
where $g$ is Riemann's integrable function and $x(t)$ is the random variable with the density $\rho(x)$. By the same method, we have
\begin{equation*}
EY=(b-a)\int_{-\infty}^{+\infty}g(x)\rho(x)\mathrm{d}x,
\end{equation*}
\begin{equation*}
DY=0,
\end{equation*}
\begin{equation*}
\rho_Y(y)=\delta(y-EY).
\end{equation*}

We know that the density of the functional  $Y=f(x)=\int_a^b g(x(t))\mathrm{d}t$ is $\rho_Y(y)=\delta(y-EY)$, so
for the general functional $Z=h(Y)$, we have formula
\begin{equation*}
EZ=Eh(Y)=h(EY),
\end{equation*}
where $h(Y)$ is a general function of $Y$. In fact,
\begin{equation*}
EZ=Eh(Y)=\int h(y)\delta(y-EY)\mathrm{d}y=h(EY).
\end{equation*}
More generally, we have the following formula.

\textbf{ Theorem 5.1}. Let $I_1,\cdots,I_k$ be $k$ subintervals of [0,1], and
$Y_1=\int_{I_1} g_1(x(t))\mathrm{d}t,\cdots,Y_k=\int_{I_k} g_k(x(t))\mathrm{d}t$,
where $g_1,\cdots,g_k$ are Riemann's integrable functions and $x(t)$ is the random variable with the density $\rho(x)$. Furthermore, Let $h(y_1,\cdots,y_k)$
be a general function. Then we have
\begin{equation*}
Eh(Y_1,\cdots,Y_k)=h(EY_1,\cdots,EY_k),
\end{equation*}
where $EY_i=\mu(I_i)\int_{-\infty}^{+\infty}g_i(x)\rho(x)\mathrm{d}x$ and $\mu(I_i)$ is the measure of $I_i$.

For example, we take $\rho(x)=\frac{1}{\sqrt{2\pi}\sigma}\exp\frac{(x-\mu)^2}{2\sigma^2}$ and $Y=\int_0^1x^2(t)\mathrm{d}t$. Then, we have
\begin{equation*}
Y_n=f_n(x)=\frac{1}{n}\sum_{k=0}^nx^2(\frac{k}{n})=\frac{1}{n}\sum_{k=0}^nx^2_k,
\end{equation*}
so
\begin{equation*}
EY_n=\frac{1}{n}\sum_{k=0}^nEx^2_k=Ex_1=\int x^2\rho(x)\mathrm{d}x=\mu^2+\sigma^2,
\end{equation*}
and
\begin{equation*}
EY=\lim_{n\rightarrow \infty}E(Y_n)=\mu^2+\sigma^2.
\end{equation*}
By
\begin{equation*}
D(Y_n)=\frac{1}{n^2}\sum_{k=0}^nD(x_k)=\frac{1}{n}D(x_1)\frac{\sigma^2}{n},
\end{equation*}
it follows that
\begin{equation*}
D(Y)=\lim_{n\rightarrow \infty}D(Y_n)=0,
\end{equation*}
and
\begin{equation*}
\rho_y(y)=\delta(y-\mu^2-\sigma^2).
\end{equation*}
For example, if $Z=\sin Y=\sin \{\int_0^1x^2(t)\mathrm{d}t\}$, we have $EZ=\sin(\mu^2+\sigma^2)$.

\chapter{The mean value of functional on subspace with codimension 1 on
$C[0,1]$}

\section{The mean value of functional
$Y=f(x)=\int_0^1a(t)x(t)\mathrm{d}t$ on the subspace
$W=\{x|\int_0^1x(t)\mathrm{d}t=s\}$ of $C[0,1]$}

 On $R[0,1]$, we consider the functional
 \begin{equation*}
Y=f(x)=\int_0^1a(t)x(t)\mathrm{d}t,
 \end{equation*}
 with the constraint
 \begin{equation*}
\int_0^1x(t)\mathrm{d}t=s.
 \end{equation*}
 By discretization, we have
\begin{equation*}
Y_n=\frac{1}{n}\sum_{k=1}^na_kx_k,
\end{equation*}
and the corresponding subspace is $M'_n=\{(x_1,\cdots,x_n)|\sum_{k=1}^nx_k=ns\}$. So
\begin{equation*}
Y_n=\sum_{k=1}^n\frac{a_k}{n}x_k,
\end{equation*}
with
\begin{equation*}
\sum_{k=1}^n x_k=ns,
\end{equation*}
and $0\leq x_k\leq ns$. Therefore, the mean value of $Y_n$ is
\begin{equation*}
EY_n=\sum_{k=1}^n\frac{a_k}{n}E(x_k)=E(x_1)\sum_{k=1}^n\frac{a_k}{n}.
\end{equation*}
We only need to compute $E(x_1)$. There are two methods to get its
value $E(x_1)=s$. One simple method is from
$E(x_1)=\cdots=E(x_n)$ and $\sum_{k=1}^n x_k=ns$. Another method is
to direct compute the corresponding integrals. Denote
$M=\{(x_1,\cdots,x_n)|\sum_{k=1}^n x_k=ns,  x_k\geq 0\}$and
$M_1=\{(x_1,\cdots,x_n)|\sum_{k=1}^{n-1} x_k\leq ns,  x_k\geq 0\}$. We
need the following results whose proofs are easy.

\textbf{Lemma 6.1}. We have
\begin{equation*}
\int_M\mathrm{d}V=\frac{\sqrt n}{(n-1)!}(ns)^{n-1},
\end{equation*}
\begin{equation*}
\int_Mx_k\mathrm{d}V=\frac{\sqrt n}{n!}(ns)^{n},
\end{equation*}
\begin{equation*}
\int_Mx_kx_j\mathrm{d}V=\frac{\sqrt n}{(n+1)!}(ns)^{n+1}, i\neq k.
\end{equation*}
\begin{equation*}
\int_Mx^2_k\mathrm{d}V=\frac{2\sqrt n}{(n+1)!}(ns)^{n+1}.
\end{equation*}
\textbf{Proof}. We only consider $x_1$ and $x_2$. By direct computation, we have
\begin{equation*}
\int_M\mathrm{d}V=\int_{M_1}\sqrt{1+(\frac{\partial z_n}{\partial z_1})^2+\cdots+(\frac{\partial z_n}{\partial z_{n-1}})^2}\mathrm{d}z_1\cdots\mathrm{d}z_{n-1}
\end{equation*}
\begin{equation*}
=\int_{M_1}\sqrt{n}\mathrm{d}z_1\cdots\mathrm{d}z_{n-1}=\frac{\sqrt n}{(n-1)!};
\end{equation*}

\begin{equation*}
\int_M z_1\mathrm{d}V=\int_{M_1}z_1\sqrt{1+(\frac{\partial z_n}{\partial z_1})^2+\cdots+(\frac{\partial z_n}{\partial z_{n-1}})^2}\mathrm{d}z_1\cdots\mathrm{d}z_{n-1}
\end{equation*}
\begin{equation*}
=\int_{M_1}\sqrt{n}z_1\mathrm{d}z_1\cdots\mathrm{d}z_{n-1}=\frac{\sqrt n}{(n-2)!}\int_0^1z_1(1-z_1)^{n-2}\mathrm{d}z_1=\frac{\sqrt n}{n!}.
\end{equation*}
\begin{equation*}
\int_Mz_1z_2\mathrm{d}V=\int_{M_1}z_1z_2\sqrt{1+(\frac{\partial z_n}{\partial z_1})^2+\cdots+(\frac{\partial z_n}{\partial z_{n-1}})^2}\mathrm{d}z_1\cdots\mathrm{d}z_{n-1}
\end{equation*}
\begin{equation*}
=\frac{\sqrt n}{(n-3)!}\int_0^1z_1\mathrm{d}z_1\int_0^{1-z_1}\mathrm{d}z_2(1-z_1-z_2)^{n-3}=\frac{\sqrt n}{(n+1)!};
\end{equation*}
and
\begin{equation*}
\int_M z^2_1\mathrm{d}V=\int_{M_1}z^2_1\sqrt{1+(\frac{\partial z_n}{\partial z_1})^2+\cdots+(\frac{\partial z_n}{\partial z_{n-1}})^2}\mathrm{d}z_1\cdots\mathrm{d}z_{n-1}
\end{equation*}
\begin{equation*}
=\int_{M_1}\sqrt{n}z^2_1\mathrm{d}z_1\cdots\mathrm{d}z_{n-1}=\frac{\sqrt n}{(n-2)!}\int_0^1z^2_1(1-z_1)^{n-2}\mathrm{d}z_1=\frac{2\sqrt n}{(n+1)!}.
\end{equation*}

By the above lemmas, we have

\textbf{Lemma 6.2}.

\begin{equation*}
E(x_k)=\frac{\sqrt n}{n!}(ns)^{n-1}/\frac{\sqrt n}{(n-1)!}(ns)^{n}=s,
\end{equation*}
\begin{equation*}
E(x_jx_k)=\frac{\sqrt n(ns)^{n+1}}{(n+1)!}/\frac{\sqrt n(ns)^{n-1}}{(n-1)!}=\frac{ns^2}{n+1}, j\neq k,
\end{equation*}
\begin{equation*}
E(x^2_k)=\frac{2\sqrt n(ns)^{n+1}}{(n+1)!}/\frac{\sqrt n(ns)^{n-1}}{(n-1)!}=\frac{2ns^2}{n+1}.
\end{equation*}

By these results, we can give $E(Y)$ and $D(Y)$.

\textbf{Theorem 6.1}.
\begin{equation*}
E(Y)=\lim_{n\rightarrow +\infty}E(Y_n)=s\int_0^1a(t)\mathrm{d}t,
\end{equation*}
\begin{equation*}
D(Y)=0.
\end{equation*}
Further, we have
\begin{equation*}
E(h(Y))=h(s\int_0^1a(t)\mathrm{d}t),
\end{equation*}
where $h(Y)$ is a general function.

\textbf{Proof.} Firstly, we have
\begin{equation*}
E(Y_n)=\frac{s}{n}\sum_{k=1}^na_k,
\end{equation*}
and
\begin{equation*}
E(Y)=\lim_{n\rightarrow +\infty}E(Y_n)=s\int_0^1a(t)\mathrm{d}t.
\end{equation*}
Now we compute the variances of $Y_n$ and $Y$. In fact, we have
\begin{equation*}
EY^2_n=\sum_{j=1}^n\sum_{k=1}^n\frac{a_ja_k}{n^2}E(x_kx_j)=\frac{1}{n^2}\sum_{k=1}^n a_k^2E(x_k^2)+\frac{1}{n^2}\sum_{k\neq j}a_ka_jE(x_kx_j)
\end{equation*}
\begin{equation*}
=\frac{1}{n^2}\sum_{k=1}^n a_k^2\frac{2ns^2}{n+1}+\frac{1}{n^2}\sum_{k\neq j}a_ka_j\frac{ns^2}{n+1}
\end{equation*}
\begin{equation*}
=\frac{s^2}{n(n+1)}\{\sum_{k=1}^n a_k\sum_{k=1}^n a_j+\sum_{k=1}^n a_k^2\},
\end{equation*}
and
\begin{equation*}
E(Y^2)=\lim_{n\rightarrow +\infty}\frac{s^2}{n(n+1)}\{\sum_{k=1}^n a_k\sum_{k=1}^n a_j+\sum_{k=1}^n a_k^2\}=(s^2(\int_0^1a(t)\mathrm{d}t)^2.
\end{equation*}
So
\begin{equation*}
D(Y)=0,
\end{equation*}
which means that the density of $Y$ is $\rho(y)=\delta(y-s\int_0^1a(t)\mathrm{d}t)$. Therefore, we have
\begin{equation*}
E(h(Y))=h(s\int_0^1a(t)\mathrm{d}t),
\end{equation*}
where $h(Y)$ is a general function. The proof is completed.

\section{The mean value of functional $Y=\int_0^1x(t)\mathrm{d}t$ on
subspace $M=\{x|\int_0^1a(t)x(t)\mathrm{d}t=s\}$ in $C[0,1]$}

By discretization, we have
\begin{equation*}
Y_n=\frac{1}{n}\sum_{k=1}^nx_k,
\end{equation*}
and the corresponding subspace is $M_n:\sum_{k=1}^na_kx_k=h$, where
$h=ns$. Therefore, the mean value of $Y_n$ is
\begin{equation*}
EY_n=\frac{1}{n}\sum_{k=1}^nE(x_k).
\end{equation*}
In order to compute $E(Y)$ and $D(Y)$, we need the following results.
Denote $M=\{(x_1,\cdots,x_n)|a_1x_1+\cdots+a_nx_n=h,a_i\geq 0,x_i\geq 0,i=1,\cdots,n\}$. Then we have
\begin{equation*}
\int_M\mathrm{d}V=\frac{\sqrt {a^2_1+\cdots+a^2_n}}{a_1\cdots a_n}\frac{h^{n-1}}{(n-1)!},
\end{equation*}
\begin{equation*}
\int_Mx_i\mathrm{d}V=\frac{\sqrt {a^2_1+\cdots+a^2_n}}{a_1\cdots a_n}\frac{1}{a_i}\frac{h^n}{n!},
\end{equation*}
\begin{equation*}
\int_Mx_ix_j\mathrm{d}V=\frac{\sqrt {a^2_1+\cdots+a^2_n}}{a_1\cdots a_n}\frac{1}{a_ia_j}\frac{h^{n+1}}{(n+1)!}, (i\neq j),
\end{equation*}
\begin{equation*}
\int_Mx^2_i\mathrm{d}V=\frac{\sqrt {a^2_1+\cdots+a^2_n}}{a_1\cdots a_n}\frac{2}{a^2_i}\frac{h^{n+1}}{(n+1)!}.
\end{equation*}

By the above formulae, we have

\textbf{Lemma 6.3}.
\begin{equation*}
Ex_i=\frac{1}{a_i}\frac{h}{n},
\end{equation*}
\begin{equation*}
E(x_ix_j)=\frac{1}{a_ia_j}\frac{h^2}{n(n+1)},
\end{equation*}
\begin{equation*}
Ex^2_i=\frac{2}{a^2_i}\frac{h^2}{n(n+1)}.
\end{equation*}
Using these formulae, we obtain
\begin{equation*}
EY_n=\frac{1}{n}\sum_{i=1}^n\frac{1}{a_i}\frac{h}{n}=\frac{s}{n}\sum_{i=1}^n\frac{1}{a_i}.
\end{equation*}
So, We have

\textbf{Theorem 6.2}.
\begin{equation*}
E(Y)=s\int_0^1\frac{1}{a(t)}\mathrm{d}t,
\end{equation*}
\begin{equation*}
D(Y)=0.
\end{equation*}
Further, we have
\begin{equation*}
E(h(Y))=h(s^2\int_0^1\frac{1}{a(t)}\mathrm{d}t),
\end{equation*}
where $h(Y)$ is a general function.

\textbf{Proof.} Using these formulae, we obtain
\begin{equation*}
EY_n=\frac{1}{n}\sum_{i=1}^n\frac{1}{a_i}\frac{h}{n}=\frac{s}{n}\sum_{i=1}^n\frac{1}{a_i}.
\end{equation*}
So,  We have
\begin{equation*}
E(Y)=s\int_0^1\frac{1}{a(t)}\mathrm{d}t.
\end{equation*}
Further, we have
\begin{equation*}
EY^2_n=\frac{1}{n^2}\sum_{j=1}^n\sum_{k=1}^nE(x^2_kx^2_j)=\frac{1}{n^2}\sum_{k=1}^n E(x_k^2)+\frac{1}{n^2}\sum_{k\neq j}E(x_kx_j)
\end{equation*}
\begin{equation*}
=\frac{1}{n^2}\sum_{k=1}^n \frac{2}{a^2_k}\frac{h^2}{n(n+1)}+\frac{1}{n^2}\sum_{k\neq j}\frac{1}{a_ka_j}\frac{h^2}{n(n+1)}
\end{equation*}
\begin{equation*}
=\sum_{k=1}^n \frac{2}{a^2_k}\frac{s^2}{n(n+1)}+\sum_{k\neq j}\frac{1}{a_ka_j}\frac{s^2}{n(n+1)},
\end{equation*}
and then
\begin{equation*}
E(Y^2)=\lim_{n\rightarrow +\infty}EY^2_n=s^2(\int_0^1\frac{1}{a(t)}\mathrm{d}t)^2,
\end{equation*}
\begin{equation*}
D(Y)=0.
\end{equation*}
This means that the density of $Y$ is $\rho(y)=\delta(y-s^2\int_0^1\frac{1}{a(t)}\mathrm{d}t)$. Therefore, we have
the nonlinear exchange formula
\begin{equation*}
E(h(Y))=h(s^2\int_0^1\frac{1}{a(t)}\mathrm{d}t),
\end{equation*}
where $h(Y)$ is a general function. The proof is completed.

\section{The mean value of functional $Y=\int_0^1g(x(t))\mathrm{d}t$ on
subspace $M=\{x|\int_0^1a(t)x(t)\mathrm{d}t=s\}$ in $C[0,1]$}

In this subsection, we use a more general method to treat a general integral form functional. Our method is to first compute the density of $x(t)$ in $M$ for the fixed $t$. By discretization, we have
\begin{equation*}
Y_n=\frac{1}{n}\sum_{k=1}^ng(x_k),
\end{equation*}
and the corresponding subspace is $M_n=\{(x_1,\cdots,x_n)|\sum_{k=1}^na_kx_k=h\}$, where
$h=ns$. Further, denote  $M'=\{(x_1,\cdots,x_n)|a_1x_1+\cdots+a_nx_n=h,a_i\geq 0,x_i\geq 0,i=1,\cdots,n\}$ and $M_1'=\{(x_2,\cdots,x_n)|a_2x_2+\cdots+a_nx_n=ns-a_1x_1,a_i\geq 0,x_i\geq 0,i=1,\cdots,n\}$. Then we have

\textbf{Lemma 6.4}. The density of $x(t)$ in $M$ for the fixed $t$ is given by
\begin{equation}
\rho(x)=\frac{a(t)}{s}\mathrm{e}^{-\frac{a(t)x}{s}}.
\end{equation}

\textbf{Proof}. We only prove the result for $x(0)$. By discretization, the density of $x_1$ is
\begin{equation*}
\rho_n(x_1)=\frac{\int_{M'_1}\mathrm{d}V}{\int_{M'}\mathrm{d}V}=\frac{\frac{\sqrt {a^2_2+\cdots+a^2_n}}{a_2\cdots a_n}\frac{(ns-a_1x_1)^{n-2}}{(n-2)!}}{\frac{\sqrt {a^2_1+\cdots+a^2_n}}{a_1\cdots a_n}\frac{(ns)^{n-1}}{(n-1)!}}
\end{equation*}
\begin{equation}
=\frac{a_1(n-1)}{sn}\frac{\sqrt {a^2_2+\cdots+a^2_n}}{\sqrt {a^2_1+\cdots+a^2_n}}(1-\frac{a_1x_1}{ns})^{n-2}.
\end{equation}
So we have
\begin{equation}
\rho(x(0))=\lim_{n\rightarrow+\infty}\frac{a_1(n-1)}{sn}\frac{\sqrt {a^2_2+\cdots+a^2_n}}{\sqrt {a^2_1+\cdots+a^2_n}}(1-\frac{a_1x_1}{ns})^{n-2}=\frac{a(0)}{s}\mathrm{e}^{-\frac{a(0)x}{s}}.
\end{equation}

\textbf{Theorem 6.3}. For the general functional $Y=\int_0^1g(x(t))\mathrm{d}t$, we have
\begin{equation*}
E(Y)=\frac{1}{s}\int_0^1\int_0^{+\infty}a(t)g(x)\mathrm{e}^{-\frac{a(t)x}{s}}\mathrm{d}x\mathrm{d}t=
\frac{1}{s}\int_0^1\int_0^{+\infty}g(\frac{sx}{a(t)})\mathrm{e}^{-x}\mathrm{d}x\mathrm{d}t,
\end{equation*}
\begin{equation*}
D(Y)=0.
\end{equation*}
Further, we have
\begin{equation*}
E(h(Y))=h(EY),
\end{equation*}
where $h(Y)$ is a general function.

\textbf{Proof.} It is easy to prove the result from the lemma 6.4. The proof is completed.

\textbf{Remark 6.1}. It is easy to see that theorems 6.1 and 6.2 are the special cases.

\chapter{The mean values of functionals on subspace with codimension 2 on
$C[0,1]$}

We will treat this problem by two methods.
\section{The first method}
Under two constraints
\begin{equation*}
\int_0^1a(t)x(t)\mathrm{d}t=r,
\end{equation*}
 \begin{equation*}
\int_0^1b(t)x(t)\mathrm{d}t=s,
\end{equation*}
we want to compute the mean value of functional $Y=\int_0^1c(t)x(t)\mathrm{d}t$, where $a(t)\geq0, b(t)\geq 0$ and $x(t)\geq0$. After discretization, we have respectively
\begin{equation*}
\sum_{i=1}^na_ix_i=nr,
\end{equation*}
 \begin{equation*}
\sum_{i=1}^nb_ix_i=ns,
\end{equation*}
and
\begin{equation*}
Y_n=\frac{1}{n}\sum_{i=1}^nc_ix_i.
\end{equation*}

Eliminating $x_n$ from two constraints and $Y_n$, these two constraints become only one
\begin{equation*}
M=\{\sum_{k=1}^{n-1}(a_nb_k-a_kb_n)x_k=n(b_nr-a_ns)\},
\end{equation*}
and corresponding $Y_n$ becomes
\begin{equation*}
Y_n=\frac{r}{a_n}+\frac{1}{na_n}\sum_{k=1}^{n-1}(a_nc_k-a_k)x_k.
\end{equation*}
In order to get $EY$ and $DY$, we need the following results:

\textbf{Lemma 7.1}.
\begin{equation*}
\int_M\mathrm{d}V=\frac{\sqrt {p^2_1+\cdots+p^2_n}}{p_1\cdots p_n}\frac{(n(b_nr-a_ns))^{n-2}}{(n-2)!},
\end{equation*}
\begin{equation*}
\int_Mx_k\mathrm{d}V=\frac{\sqrt {p^2_1+\cdots+p^2_n}}{p_1\cdots p_n}\frac{(n(b_nr-a_ns))^{n-1}}{(a_kb_n-a_nb_k)(n-1)!},
\end{equation*}
\begin{equation*}
\int_Mx_kx_j\mathrm{d}V=\frac{\sqrt {p^2_1+\cdots+p^2_n}}{p_1\cdots p_n}\frac{(n(b_nr-a_ns))^{n+1}}{(a_kb_n-a_nb_k)(a_jb_n-a_nb_j)(n+1)!}, (i\neq j),
\end{equation*}
\begin{equation*}
\int_Mx^2_k\mathrm{d}V=\frac{\sqrt {p^2_1+\cdots+p^2_n}}{p_1\cdots p_n}\frac{2(n(b_nr-a_ns))^{n+1}}{(a_kb_n-a_nb_k)(n+1)!},
\end{equation*}
where $p_k=a_kb_n-a_nb_k$.

By the above formulae, we have

\textbf{Lemma 7.2}.
\begin{equation*}
Ex_k=\frac{b_nr-a_ns}{a_kb_n-a_nb_k}\frac{n}{n-1},
\end{equation*}
\begin{equation*}
E(x_kx_j)=\frac{(b_nr-a_ns)^2}{a_kb_n-a_nb_k}\frac{2n}{n+1},
\end{equation*}
\begin{equation*}
Ex^2_k=\frac{(b_nr-a_ns)^2}{(a_kb_n-a_nb_k)(a_jb_n-a_nb_j)}\frac{n}{n+1}.
\end{equation*}

\textbf{Theorem 7.1.}
\begin{equation*}
E(Y)=\frac{r}{a(1)}+\frac{b(1)r-a(1)s}{a(1)}\int_0^1\frac{a(1)c(t)-a(t)}{a(1)b(t)-a(t)b(1)}\mathrm{d}t,
\end{equation*}
\begin{equation*}
D(Y)=0.
\end{equation*}
Further,  we have the nonlinear exchange formula
\begin{equation*}
E(h(Y))=h(\frac{r}{a(1)}+\frac{(b(1)r-a(1)s)}{a(1)}\int_0^1\frac{a(1)c(t)-a(t)}{a(1)b(t)-a(t)b(1)}\mathrm{d}t),
\end{equation*}
where $h(Y)$ is a general function.

\textbf{Proof.} By lemma 7.2, we obtain
\begin{equation*}
EY_n=\frac{r}{a_n}+\frac{b_nr-a_ns}{a_n}\frac{1}{n-1}\sum_{k=1}^{n-1}\frac{a_nc_k-a_k}{a_nb_k-a_kb_n},
\end{equation*}
ad then,
\begin{equation*}
E(Y)=\frac{r}{a(1)}+\frac{b(1)r-a(1)s}{a(1)}\int_0^1\frac{a(1)c(t)-a(t)}{a(1)b(t)-a(t)b(1)}\mathrm{d}t.
\end{equation*}
Furthermore, we have
\begin{equation*}
EY^2_n=\frac{r^2}{a^2_n}+\frac{2r(b_nr-a_ns)}{na^2_n}\sum_{k=1}^{n-1}\frac{a_nc_k-a_k}{a_nb_k-a_kb_n}
\end{equation*}
\begin{equation*}
+\frac{2(b_nr-a_ns)^2}{n(n+1)a^2_n}\sum_{k=1}^{n-1}\frac{(a_nc_k-a_k)^2}{a_nb_k-a_kb_n}+\frac{2(b_nr-a_ns)^2}{n(n+1)a^2_n}\sum_{k\neq j}\frac{(a_nc_k-a_k)(a_nc_j-a_j)}{(b_na_k-a_nb_k)(b_na_j-a_nb_j)}
\end{equation*}
\begin{equation*}
=\frac{r^2}{a^2_n}+\frac{2r(b_nr-a_ns)}{na^2_n}\sum_{k=1}^{n-1}\frac{a_nc_k-a_k}{a_nb_k-a_kb_n}
\end{equation*}
\begin{equation*}
+\frac{2(b_nr-a_ns)^2}{(n+1)^2a^2_n}\{\sum_{k=1}^{n-1}\frac{a_nc_k-a_k}{a_nb_k-a_kb_n}\}^2
+\frac{(b_nr-a_ns)^2}{n(n+1)a^2_n}\sum_{k=0}^{n-1}\frac{(a_nc_k-a_k)^2}{(a_nb_k-a_kb_n)^2},
\end{equation*}
and so
\begin{equation*}
EY^2=\frac{r^2}{a^2(1)}+\frac{2r(b(1)r-a(1)s)}{a^2(1)}\int_0^1\frac{a(1)c(t)-a(t)}{a(1)b(t)-a(t)b(1)}\mathrm{d}t
\end{equation*}
\begin{equation*}
+\frac{(b(1)r-a(1)s)^2}{a^2(1)}\{\int_0^1\frac{a(1)c(t)-a(t)}{a(1)b(t)-a(t)b(1)}\mathrm{d}t\}^2
\end{equation*}
\begin{equation*}
=\{\frac{r}{a(1)}+\frac{(b(1)r-a(1)s)}{a(1)}\int_0^1\frac{a(1)c(t)-a(t)}{a(1)b(t)-a(t)b(1)}\mathrm{d}t\}^2=E^2Y,
\end{equation*}
and
\begin{equation*}
D(Y)=0.
\end{equation*}
This means that the density of $Y$ is
\begin{equation*}
\rho(y)=\delta(y-\frac{r}{a(1)}-\frac{(b(1)r-a(1)s)}{a(1)}\int_0^1\frac{a(1)c(t)-a(t)}{a(1)b(t)-a(t)b(1)}\mathrm{d}t)
\end{equation*}
Therefore, we have the nonlinear exchange formula
\begin{equation*}
E(h(Y))=h(\frac{r}{a(1)}+\frac{(b(1)r-a(1)s)}{a(1)}\int_0^1\frac{a(1)c(t)-a(t)}{a(1)b(t)-a(t)b(1)}\mathrm{d}t),
\end{equation*}
where $h(Y)$ is a general function. The proof is completed.

\section{The second method}

We consider the general functional $Y=\int_0^1g(x(t))\mathrm{d}t$. Our method is to compute the density of $x(t)$ in $M$. We denote $M'=\{(x_1,\cdots,x_{n-1})|\sum_{k=1}^{n-1}(a_nb_k-a_kb_n)x_k=n(b_nr-a_ns)\}$, and $M_1'=\{(x_2,\cdots,x_{n-1})|\sum_{k=1}^{n-1}(a_nb_k-a_kb_n)x_k=n(b_nr-a_ns)-a_nb_1+a_1b_n\}$. Then we have

\textbf{Lemma 7.3}. The density of $x(t)$ in $M$ for the fixed $t$ is given by
\begin{equation}
\rho(x)=\frac{a(t)b(1)-a(1)b(t)}{b(1)r-a(1)s}\mathrm{e}^{-\frac{a(1)b(t)-b(1)a(t)}{b(1)r-a(1)s}x(t)}.
\end{equation}

\textbf{Proof}. We only prove the result for $x(0)$. By discretization, the density of $x_1$ is
\begin{equation*}
\rho_n(x_1)=\frac{\int_{M'_1}\mathrm{d}V}{\int_{M'}\mathrm{d}V}=\frac{\frac{\sqrt {p^2_2+\cdots+p^2_n}}{p_2\cdots p_n}\frac{(n(b_nr-a_ns)-(a_nb_1-a_1b_n)x_1)^{n-3}}{(n-3)!}}{\frac{\sqrt {p^2_1+\cdots+p^2_n}}{p_1\cdots p_n}\frac{(n(b_nr-a_ns))^{n-2}}{(n-2)!}}
\end{equation*}
\begin{equation}
=\frac{(a_1b_n-a_nb_1)(n-1)}{n(b_nr-a_ns)}\frac{\sqrt {p^2_2+\cdots+p^2_n}}{\sqrt {p^2_1+\cdots+p^2_n}}(1-\frac{(b_nb_1-b_na_1)x_1}{n(b_nr-a_ns)})^{n-3}.
\end{equation}
So we have
\begin{equation}
\rho(x(0))=\lim_{n\rightarrow+\infty}\rho_n(x_1)
=\frac{a(0)b(1)-a(1)b(0)}{b(1)r-a(1)s}\mathrm{e}^{-\frac{a(1)b(0)-b(1)a(0)}{b(1)r-a(1)s}x(0)}.
\end{equation}
By replacing $0$ by $t$ in the above formula, we give the result. The proof is completed.

\textbf{Theorem 7.2}. For the general functional $Y=\int_0^1g(x(t))\mathrm{d}t$, we have
\begin{equation*}
E(Y)=\int_0^1\int_0^{+\infty}g(x)
\frac{a(t)b(1)-a(1)b(t)}{b(1)r-a(1)s}\mathrm{e}^{-\frac{a(1)b(t)-b(1)a(t)}{b(1)r-a(1)s}x(0)}\mathrm{d}x\mathrm{d}t,
\end{equation*}
\begin{equation*}
D(Y)=0.
\end{equation*}
Further, we have
\begin{equation*}
E(h(Y))=h(EY),
\end{equation*}
where $h(Y)$ is a general function.

\textbf{Proof.} It is easy to prove the result from the lemma 7.3. The proof is completed.

\textbf{Remark 7.1}. It is easy to see that theorems 7.1 is the special case.

\chapter{The mean values of functionals on the infinite-dimensional balls in $C[0,1]$ with 2-norm}

Consider the ball $M$ in function space  $C[0,1]$ with 2-norm,
\begin{equation}
M=\{x|\int_0^1x^2(t)\mathrm{d}t\leq R^2\}.
\end{equation}
We will firstly compute the mean values of the following functionals on $M$
\begin{equation}
Y=f(x)=\int_0^1x(t)\mathrm{d}t,
\end{equation}
\begin{equation}
Y=f(x)=\int_0^1x^2(t)\mathrm{d}t.
\end{equation}
Then we give the mean value of  the general functional with integral form
\begin{equation}
Y=f(x)=\int_0^1g(x(t))\mathrm{d}t.
\end{equation}
We discretize it and $M$ as follows
\begin{equation}
Y_n=\frac{1}{n}\sum_{k=1}^ng(x_k),
\end{equation}
\begin{equation}
M_n=\{(x_1,\cdots,x_n)|\sum_{k=1}^nx_k^2\leq nR^2\},
\end{equation}
where $x_k=\frac{k}{n}$. Then we define the mean value of $Y$ as
\begin{equation}
EY=\lim_{n\rightarrow\infty}E(Y_n)=\lim_{n\rightarrow\infty}
\frac{\frac{1}{n}\sum_{k=1}^n\int_{M_n}g(x_k)\mathrm{d}V}{\int_{M_n}\mathrm{d}V},
\end{equation}
where $\mathrm{d}V$ is the volume element of $M_n$. To compute the mean of values, we need the following result.

\textbf{Lemma 8.1}([39]). The following generalized Dirichlet formula holds
\begin{equation}
\int\cdots\int_{B^+}x_1^{p_1-1}x_2^{p_2-1}\cdots x_n^{p_n-1}\mathrm{d}x_1\cdots\mathrm{d}x_n=
\frac{1}{2^n}\frac{\Gamma(\frac{p_1}{2})\cdots \Gamma(\frac{p_n}{2})}{\Gamma(1+\frac{p_1+\cdots p_n}{2})},
\end{equation}
where $B^+=\{(x_1,\cdots,x_n)|x_1^2+\cdots+x_n^2\leq1, x_k>0, k=1,\cdots, n\}$.

By the above lemma, we have the following results.

\textbf{Theorem 8.1}. For the functional $Y=f(x)=\int_0^1x(t)\mathrm{d}t$, its mean value and variance on $M$ are given by
\begin{equation}
EY=0,
\end{equation}
\begin{equation}
DY=0,
\end{equation}
and hence
\begin{equation}
Eh(Y)=h(EY)=h(0),
\end{equation}
where $h$ is a general function.

\textbf{Proof}. It is easy to see that
\begin{equation}
\int\cdots\int_{B}x_k\mathrm{d}x_1\cdots\mathrm{d}x_n=0,
\end{equation}
where $B=\{(x_1,\cdots,x_n)|x_1^2+\cdots+x_n^2\leq1\}$. So $EY=0$. Then we prove $DY=0$. Indeed, by symmetry of $B$, we have for $i\neq j$,
\begin{equation}
\int\cdots\int_{B}x_ix_j\mathrm{d}x_1\cdots\mathrm{d}x_n=0.
\end{equation}
Furthermore, by the lemma 8.1, we have
\begin{equation}
E(x_k^2)=\frac{\int\cdots\int_{M_n}x_k^2\mathrm{d}x_1\cdots\mathrm{d}x_n}{\int_{M_n}\mathrm{d}V}
=\frac{nR^2}{2+n}.
\end{equation}
Therefore, we have
\begin{equation}
E(Y^2_n)=\frac{1}{n^2}\{\sum_{k=1}^nE(x^2_k)+\sum_{i\neq j}E(x_ix_j)\}=\frac{nR^2}{n(2+n)},
\end{equation}
from which we get $E(Y^2)=0$, and hence $DY=0$. The last conclusion is obvious. The proof is completed.

\textbf{Theorem 8.2}. For the functional $Y=f(x)=\int_0^1x^2(t)\mathrm{d}t$, its mean value and variance on $M$ are given by
\begin{equation}
EY=R^2,
\end{equation}
\begin{equation}
DY=0,
\end{equation}
and hence
\begin{equation}
Eh(Y)=h(EY)=h(R^2),
\end{equation}
where $h$ is a general function.

\textbf{Proof}. By the lemma 8.1, we have
\begin{equation}
E(x_k^2)=\frac{\frac{\frac{1}{2^n}n^{\frac{n}{2}}R^n(\sqrt nR)^2\Gamma(\frac{3}{2})\Gamma^{n-1}(\frac{1}{2})}{\Gamma(2+\frac{n}{2})}}{\frac{n^{\frac{n}{2}}R^n\pi^{\frac{n}{2}}}
{2^n\Gamma(1+\frac{n}{2})}}=\frac{nR^2}{n+2},
\end{equation}
\begin{equation}
E(x_k^4)=\frac{\frac{\frac{1}{2^n}n^{\frac{n}{2}}R^n(\sqrt nR)^4\Gamma(\frac{5}{2})\Gamma^{n-1}(\frac{1}{2})}{\Gamma(3+\frac{n}{2})}}{\frac{n^{\frac{n}{2}}R^n\pi^{\frac{n}{2}}}
{2^n\Gamma(1+\frac{n}{2})}}=\frac{3n^2R^4}{(n+2)(n+4)}.
\end{equation}
and for $i\neq j$,
\begin{equation}
E(x_i^2x_j^2)=\frac{\frac{\frac{1}{2^n}n^{\frac{n}{2}}R^n(\sqrt nR)^4\Gamma^2(\frac{3}{2})\Gamma^{n-2}(\frac{1}{2})}{\Gamma(3+\frac{n}{2})}}{\frac{n^{\frac{n}{2}}R^n\pi^{\frac{n}{2}}}
{2^n\Gamma(1+\frac{n}{2})}}=\frac{n^2R^4}{(n+2)(n+4)}.
\end{equation}
Therefore, we have
\begin{equation}
E(Y_n)=\frac{1}{n}\sum_{k=1}^nE(x^2_k)=\frac{nR^2}{2+n},
\end{equation}
from which we get $E(Y)=R^2$.  Furthermore, we have
\begin{equation}
E(Y^2_n)=\frac{1}{n^2}\{\sum_{k=1}^nE(x^4_k)+\sum_{i\neq j}E(x^2_ix^2_j)\}=\frac{3nR^4}{(n+2)(n+4)}+\frac{n(n-1)R^4}{(n+2)(n+4)},
\end{equation}
from which we get $EY^2=R^4$ and hence $DY=0$. The last conclusion is obvious. The proof is completed.

In general, we have the following theorem.

\textbf{Theorem 8.3}. For the functional $Y=f(x)=\int_0^1g(x(t))\mathrm{d}t$, its mean value and variance on $M$ are given by
\begin{equation}
EY=\frac{1}{\sqrt{2\pi}}\int_{-\infty}^{+\infty}g(x)\mathrm{e}^{-\frac{x^2}{2}}\mathrm{d}x,
\end{equation}
\begin{equation}
DY=0,
\end{equation}
and hence
\begin{equation}
Eh(Y)=h(EY),
\end{equation}
where $h$ is a general function.

\textbf{Remark 8.1}. We will give two proofs to this theorem. The first proof is based on pure computation under the assumption of $g(x)$ being an analytic function. The second proof is based on more general idea and offer an insight to the essence of this kind of problems and hence can be generalize to other similar problems which will be considered in next section.

\textbf{The first proof}. It is easy to see that we have on $M_n$,
\begin{equation}
\int\cdots\int_{B}x^{2m+1}_k\mathrm{d}x_1\cdots\mathrm{d}x_n=0,
\end{equation}
and by the lemma 8.1,
\begin{equation}
E(x_k^{2m})=\frac{\int\cdots\int_{M_n}x_k^{2m}\mathrm{d}x_1\cdots\mathrm{d}x_n}{\int_{M_n}\mathrm{d}V}
=\frac{n^mR^{2m}\Gamma(m+\frac{1}{2}\Gamma(1+\frac{n}{2}))}{\sqrt{\pi}\Gamma(m+1+\frac{n}{2})}.
\end{equation}
Therefore, on $M$, we have
\begin{equation}
E(x_k^{2m})=\lim_{n\rightarrow+\infty}
=\frac{n^mR^{2m}(2m-1)!!}{(n+2m)(n+2m-2)\cdots(n+2)}=(2m-1)!!R^{2m}.
\end{equation}
Then, by assumption of $g(x)$ being an analytic function, we have on $M$
\begin{equation}
g(x_k)=\sum_{m=0}^{+\infty}\frac{g^{(m)}(0)}{m!}x_k^m,
\end{equation}
and hence
\begin{equation*}
E(g(x_k))=\sum_{m=0}^{+\infty}\frac{g^{(m)}(0)}{m!}E(x_k^m)=\sum_{m=0}^{+\infty}\frac{g^{(2m)}(0)}{(2m)!}E(x_k^{2m}),
\end{equation*}
\begin{equation}
=\sum_{m=0}^{+\infty}\frac{g^{(2m)}(0)}{(2m)!}(2m-1)!!R^{2m}
=\sum_{m=0}^{+\infty}\frac{g^{(2m)}(0)}{2^mm!}R^{2m}.
\end{equation}
On the other hand, we have
\begin{equation*}
\int_{-\infty}^{+\infty}g(Rx)\mathrm{e}^{-\frac{x^2}{2}}\mathrm{d}x
=\sum_{m=0}^{+\infty}\int_{-\infty}^{+\infty}\frac{g^{(m)}(0)}{m!}R^mx^m\mathrm{e}^{-\frac{x^2}{2}}\mathrm{d}x
\end{equation*}
\begin{equation}
=\sum_{m=0}^{+\infty}\frac{g^{(2m)}(0)R^{2m}}{(2m)!}\int_{-\infty}^{+\infty}x^{2m}\mathrm{e}^{-\frac{x^2}{2}}\mathrm{d}x
=\sqrt{2\pi}\sum_{m=0}^{+\infty}\frac{g^{(2m)}(0)}{(2m)!}(2m-1)!!R^{2m}.
\end{equation}
It follows that on $M$
\begin{equation}
E(g(x_k))=\frac{1}{\sqrt{2\pi}}\int_{-\infty}^{+\infty}g(Rx)\mathrm{e}^{-\frac{x^2}{2}}\mathrm{d}x.
\end{equation}
Therefore, by the linearity of the expectation $E$,
\begin{equation}
EY=\int_0^1E(g(x(t))\mathrm{d}t=E(g(x(t))=E(g(x_k))
=\frac{1}{\sqrt{2\pi}}\int_{-\infty}^{+\infty}g(Rx)\mathrm{e}^{-\frac{x^2}{2}}\mathrm{d}x.
\end{equation}
Similarly, we have,
\begin{equation}
E(Y^2)=\int_0^1\int_0^1E(g(x(t))E(g(x(s))\mathrm{d}t\mathrm{d}s=E(g(x(t)g(x(s))=E(g(x_kx_j)).
\end{equation}
By further computation, we can obtain
\begin{equation}
E(Y^2)=\int_0^1\int_0^1E(g(x(t))E(g(x(s))\mathrm{d}t\mathrm{d}s=E(g(x(t))E(g(x(s))=E^2(Y),
\end{equation}
and then $DY=0$. The proof is completed.

\textbf{The second proof}. We direct compute the density $\rho_n(x)$ of $x_1$ as a random variable in the ball $M_n$. We have
\begin{equation*}
\rho_n(x_1)=\frac{\int\cdots\int_{M_n}\mathrm{d}x_2\cdots\mathrm{d}x_n}{\int_{M_n}\mathrm{d}V}
=\frac{\frac{\pi^{\frac{n-1}{2}}}{\Gamma(1+\frac{n-1}{2})}(nR^2-x^2_1)^{\frac{n-1}{2}}}
{\frac{\pi^{\frac{n}{2}}}{\Gamma(1+\frac{n}{2})}(nR^2)^{\frac{n}{2}}}
\end{equation*}
\begin{equation}
=\frac{\Gamma(1+\frac{n}{2})}{\sqrt\pi\Gamma(1+\frac{n-1}{2})\sqrt{nR^2-x^2_1}}(1-\frac{x_1^2}{nR^2})^{\frac{n}{2}}.
\end{equation}
Taking the limitation of $n$ approaching to $+\infty$, and using the Stirling's asymptotic formula of Gamma function,  we get for $t$ fixed the density of $x(t)$ in $M$,
\begin{equation*}
\rho(x)=\lim_{n\rightarrow+\infty}\frac{\sqrt{2\pi}\sqrt{\frac{n}{2}}
(\frac{n}{2})^{\frac{n}{2}}\mathrm{e}^{-\frac{n}{2}}}
{\sqrt\pi\sqrt{2\pi}\sqrt{\frac{n-1}{2}}
(\frac{n-1}{2})^{\frac{n-1}{2}}\mathrm{e}^{-\frac{n-1}{2}}\sqrt{nR^2-x^2}}(1-\frac{x^2}{nR^2})^{\frac{n}{2}}
\end{equation*}
\begin{equation}
=\lim_{n\rightarrow+\infty}\frac{1}{\sqrt \mathrm{e}\sqrt{2\pi}}\sqrt{\frac{n}{nR^2-x^2}}(\frac{n}{n-1})^{\frac{n}{2}}
(1-\frac{x^2}{nR^2})^{\frac{n}{2}}=\frac{1}{R\sqrt{2\pi}}\mathrm{e}^{-\frac{x^2}{2R^2}}.
\end{equation}
If $R=1$, it follows that in the unit ball of $R^2[0,1]$ the density of $x(t)$ as a random variable for $t$ fixed is exactly the standard normal distribution. So we have
\begin{equation}
E(g(x(t)))=\frac{1}{R\sqrt{2\pi}}\int_{-\infty}^{+\infty}g(x)\mathrm{e}^{-\frac{x^2}{2R^2}}\mathrm{d}x
=\frac{1}{\sqrt{2\pi}}\int_{-\infty}^{+\infty}g(Rx)\mathrm{e}^{-\frac{x^2}{2}}\mathrm{d}x.
\end{equation}
Other discussions are similar with the first proof. The proof is completed.

For the functional represented by a multiple integration, we can easily prove the following theorem.

\textbf{Theorem 8.4}. Take the ball $M=\{x|\int_0^1x^p(t)\mathrm{d}t\leq R^p\}$ in $C[0,1]$ with norm $||x||_p$ when  $p=\frac{p_0}{q_0}$ where $p_0$ is even and $(p_0,q_0)=1$. Then the mean-value and invariance of the functional
\begin{equation}
Y=f(x)=\int_0^1\cdots \int_0^1g(x(t_1),\cdots, x(t_m))\mathrm{d}t_1\cdots\mathrm{d}t_m,
\end{equation}
satisfy
\begin{equation}
EY=\frac{1}{(2\pi)^{\frac{m}{2}}}\int_{-\infty}^{+\infty}\cdots\int_{-\infty}^{+\infty}
g(Rx_1,\cdots,Rx_m)\mathrm{e}^{-\frac{x_1^2}{2}-\cdots-\frac{x_m^2}{2}}\mathrm{d}x_1\cdots\mathrm{d}x_m,
\end{equation}
\begin{equation}
DY=0.
\end{equation}
Further, for the functionals $Y_1,\cdots, Y_m$ with the form (8.40) and a general function $h$ of $m$ variables,  we have the nonlinear exchange formula for the mean-value of $h(Y_1,\cdots,Y_m)$ on $M$,
\begin{equation}
Eh(Y_1,\cdots,Y_m)=h(EY_1,\cdots,EY_m).
\end{equation}

\textbf{Remark 8.2}. For the ball in space $C[0,1]$ with $p-$norm, the probability density of every coordinate of any point in the ball has been given in my paper{40].

\chapter{Cauchy space}

\section{Cauchy space and its properties}
For an integral form functional $Y=f(x)$, in general, its variance
is zero, i.e.,$DY=0$ if the second moment of $x(t)$ as a random
variable is finite. This means that the measure of the set
$\{x(t)|f(x)=EY\}$ is 1. In other words, the measure of the set
$E_{ab}=\{a<f(x)\leq b\}$ is 1 or 0 according to $EY\in(a,b]$ or
$EY\bar\in(a,b]$. We can not obtain a set $E_{ab}$ such that
$0<\mu(E_{ab})<1$ for the functional $Y$. We call such set $E_{ab}$ with $0<\mu(E_{ab})<1$
the nontrivial measurable set. In order to find the nontrivial
measurable set, we introduce the following Cauchy space.

\textbf{Definition 9.1}. As a random variable, if the density of $x(t)$ is Cauchy distribution
\begin{equation*}
\rho(x)=\frac{1}{\pi}\frac{1}{1+x^2},
\end{equation*}
then the space of Riemman's integrable function $R[0,1]$ is called the Cauchy space.

We know that for Cauchy random variable $x(t)$, its characteristic function is $\phi(t)=\mathrm{e}^{|t|}$.
Hence, for the sum of Cauchy random variables $x_k$,
\begin{equation*}
Y_n=\frac{1}{n}\sum_{k=1}^n x_k,
\end{equation*}
its characteristic function is
\begin{equation*}
\phi_n(t)=\phi^n(\frac{t}{n})=(\mathrm{e}^{|\frac{t}{n}|})^n=\mathrm{e}^{|t|}.
\end{equation*}
It follows that the characteristic function of the functional
\begin{equation*}
Y=f(x)=\int_0^1x(t)\mathrm{d}t
\end{equation*}
is also
\begin{equation*}
\phi(t)=\mathrm{e}^{|t|}.
\end{equation*}
This means that the density of $Y$ is also the Cauchy density
\begin{equation*}
\rho_Y(y)=\frac{1}{\pi}\frac{1}{1+y^2}.
\end{equation*}
Now it is easy to construct nontrivial measurable sets by using the functional $Y=f(x)$. For example, the measure of the set $E_{ab}$ is
\begin{equation*}
\mu(E_{ab})=\int_a^b\frac{1}{\pi}\frac{1}{1+y^2}\mathrm{d}y.
\end{equation*}

\textbf{Example 9.1}. For $Y=f(x)=\{1+(\int_0^1x(t)\mathrm{d}t)^2\}\mathrm{e}^{-(\int_0^1x(t)\mathrm{d}t)^2}$,
we have
\begin{equation*}
EY=\int_{-\infty}^{+\infty}(1+y^2)\mathrm{e}^{-y^2}\frac{1}{\pi}\frac{1}{1+y^2}\mathrm{d}y
=\frac{1}{\pi}\int_{-\infty}^{+\infty}\mathrm{e}^{-y^2}\mathrm{d}y=\frac{1}{\sqrt{\pi}},
\end{equation*}
\begin{equation*}
EY^2=\frac{1}{\pi}\int_{-\infty}^{+\infty}(1+y^2)\mathrm{e}^{-2y^2}\mathrm{d}y=(\frac{\sqrt2}{8}+\frac{1}{2})\frac{1}{\sqrt{\pi}},
\end{equation*}
and
\begin{equation*}
DY=(\frac{\sqrt2}{8}+\frac{1}{2})\frac{1}{\sqrt{\pi}}-\frac{1}{\pi}.
\end{equation*}

In the Cauchy space there is a rich measure theory for the functionals
of integral form. In some degree, the Cauchy space is unique. We
suppose that $x_1,\cdots,x_k,\cdots$ are the random variables with
independent identity distribution, and their characteristic function
is $\phi(t)$. For $Y=\frac{1}{n}\sum_{k=1}^n x_k$, its
characteristic function is $\phi^n(\frac{t}{n})$. If real function
$\phi(t)$ is not constant 1, and
\begin{equation*}
\phi(t)=\phi^n(\frac{t}{n})
\end{equation*}
then we can almost conclude that
\begin{equation*}
\phi(t)=\mathrm{e}^{|t|}.
\end{equation*}

In fact, we only consider the case of  $n=2$, and the general case can be dealt with by similar method.
By $\phi(t)=\phi^2(\frac{t}{2})$, we have $\phi(2t)=\phi^2(t)$. Furthermore, we assume that for $t\geq 0$,
 $\phi(t)$ can be expanded as a series
 \begin{equation*}
 \phi(t)=a_0+a_1t+a_2t^2+\cdots+a_kt^k+\cdots.
 \end{equation*}
From $\phi(0)=1$, we get $a_0=1$. Substituting the above series into the equation of $\phi$ yields
\begin{equation*}
2^ka_k=\sum_{i=0}^ka_ma_{k-m},
\end{equation*}
which follows that
\begin{equation*}
a_k=\frac{a_1}{k!},
\end{equation*}
and then
\begin{equation*}
 \phi(t)=1+a_1t+\frac{a_1}{2!}t^2+\cdots+\frac{a_1}{k!}t^k+\cdots=\mathrm{e}^{a_1t},
 \end{equation*}
where $a_1\neq 0$. Since $\phi(t)$ is real function, then $\phi(-t)=\phi(t)$,
which  gives $a_1t=|a_1t|$. We take $a_1=1$ to get the characteristic function of the Cauchy distribution.

\section{The transformation between the uniform random variable and the
Cauchy random variable}

Let $X$ be a uniform random variable on $[0,1]$. We want to find a
function $Y=g(X)$ such that $Y$ is a Cauchy random variable on $R$.
Since the density of $Y$ satisfies
\begin{equation*}
P\{x|y<g(x)\leq y+dy\}=\frac{1}{\pi}\frac{1}{1+y^2}dy,
\end{equation*}
we take $g(x)$ to be a monotone increasing function to give
\begin{equation*}
x=g^{-1}(y)=\int\frac{1}{\pi}\frac{1}{1+y^2}dy=\frac{1}{\pi}\arctan y,
\end{equation*}
that is
\begin{equation*}
Y=\tan(\pi X).
\end{equation*}
This is just the relation between the uniform variable $X$ and the Cauchy variable $Y$. Therefore, using this transformation, we can construct the nontrivial measurable set in $R[0,1]$.

\textbf{Example 9.2}. Consider $Y=\int_0^1\tan(\pi x(t))\mathrm{d}t$, where $x$ is the uniform random variable on $[0,1]$. So
$\tan(\pi x(t))$ is a Cauchy variable, and then the density of $Y$ is $\frac{1}{\pi}\frac{1}{1+y^2}$. Therefore, we know that using $Z=(1+Y^2)\mathrm{e}^{-Y^2}$ can give nontrivial measurable sets.

Furthermore, we give the corresponding concentration of measure and nonlinear exchange formula in the Cauchy space. Let
\begin{equation*}
Y=\int_{-\infty}^{+\infty}\frac{1}{\pi}\arctan(x(t))\mathrm{d}t,
\end{equation*}
where $x$ is a Cauchy random variable. Then,we have
\begin{equation*}
EY=\frac{1}{2},
\end{equation*}
\begin{equation*}
 DY=0.
\end{equation*}
Thus, for general function $h(Y)$, we obtain the nonlinear exchange formula
\begin{equation*}
Eh(Y)=h(EY).
\end{equation*}

\chapter{Some discussions on the mean values of functionals on Wiener
space}

\section{An example of the functional with nonzero variance}
 In previous sections, we discuss the mean values of functionals on some function spaces and find the
concentration of measure phenomenon. However, for the most important function space-Wiener space, the variance of the general functional such as
\begin{equation*}
Y=\int_0^1 x^2(t)\mathrm{d}t,
\end{equation*}
is not zero. This means that we can not to use the concentration of measure phenomenon to compute the mean value of more complicate functional such as
\begin{equation*}
h(Y)=\mathrm{e}^{\mathrm{i}\xi\int_0^1 x^2(t)\mathrm{d}t}.
\end{equation*}
Wiener[14-17], Cameron and Martin[36-38], Kac[26-29] had been computed the mean values of some functionals on Wiener space by some mathematical tricks such as transformation and the Sturm-Liouvlle theory.

The reason for nonexistence of the concentration of measure phenomenon is that two random variables at two time points are not independent each other. Now, we give the concrete computations of the mean value and the variance of the functional $Y$. By discretization,  we have
\begin{equation*}
Y_n=\frac{1}{n}\sum_{k=1}^n x^2(\frac{k}{n})=\frac{1}{n}\sum_{k=1}^n x^2_k.
\end{equation*}
By direct computation, we get
\begin{equation*}
E(x_k)=(\frac{n}{2\pi})^{\frac{n}{2}}\int_{-\infty}^{+\infty}\cdots\int_{-\infty}^{+\infty}x^2_k
\exp\{-\frac{n}{2}\sum_{m=1}^n(x_m-x_{m-1})^2\}\mathrm{d}x_1\cdots\mathrm{d}x_n=\frac{k}{n}.
\end{equation*}
So,
\begin{equation*}
EY_n=\frac{1}{n}\sum_{k=1}^n\frac{k}{n}=\frac{n+1}{2n},
\end{equation*}
and then
\begin{equation*}
EY=\frac{1}{2}.
\end{equation*}
In order to compute the invariance of $Y$, we first give the following results. For $k>j$,
\begin{equation*}
E(x^2_kx^2_j)=(\frac{n}{2\pi})^{\frac{n}{2}}\int_{-\infty}^{+\infty}\cdots\int_{-\infty}^{+\infty}x^2_kx^2_j
\exp\{-\frac{n}{2}\sum_{m=1}^n(x_m-x_{m-1})^2\}\mathrm{d}x_1\cdots\mathrm{d}x_n
\end{equation*}
\begin{equation*}
=\frac{kj+2j^2}{n^2}.
\end{equation*}
\begin{equation*}
E(x^4_k)=(\frac{n}{2\pi})^{\frac{n}{2}}\int_{-\infty}^{+\infty}\cdots\int_{-\infty}^{+\infty}x^4_k
\exp\{-\frac{n}{2}\sum_{m=1}^n(x_m-x_{m-1})^2\}\mathrm{d}x_1\cdots\mathrm{d}x_n=\frac{3k^2}{n^2}.
\end{equation*}
Hence, the variance of $Y_n$ is
\begin{equation*}
EY^2_n=\frac{1}{n^2}\{2\sum_{k>j}\frac{kj+2j^2}{n^2}+\sum_{k=1}^n\frac{3k^2}{n^2}\}
\end{equation*}
\begin{equation*}
=\frac{2}{n^4}\sum_{k>j}(kj+2j^2)+\frac{3}{n^4}\sum_{k=1}^nk^2
\end{equation*}
\begin{equation*}
=\frac{2}{n^4}\sum_{k=2}^nk\sum_{j=1}^{k-1}j+\frac{4}{n^4}\sum_{k=2}^n\sum_{j=1}^{k-1}j^2+\frac{3}{n^4}\sum_{k=1}^nk^2
\end{equation*}
\begin{equation*}
=\frac{7}{3n^4}\sum_{k=1}^nk^3+\cdots=\frac{7}{3n^4}\frac{n^2(n+1)^2}{4}+\cdots.
\end{equation*}
Furthermore, we have
\begin{equation*}
EY^2=\lim_{n\rightarrow \infty}EY^2_n=\frac{7}{12},
\end{equation*}
it follows that
\begin{equation*}
DY=EY^2-E^2Y=\frac{7}{12}-\frac{1}{4}=\frac{1}{3}.
\end{equation*}
This means that the nonlinear exchange formula is not right again.

\section{Transformation between dependence and independence}

If we take a transformation to transfer the dependent variables into independent variables, the reason why the variance is not zero will be seen clearly. In fact, we take
\begin{equation*}
y_k=x_k-x_{k-1}, k=1,\cdots,n.
\end{equation*}
Correspondingly, we get
\begin{equation*}
x_k=y_1+\cdots+y_k, k=1,\cdots,n.
\end{equation*}
Then
\begin{equation*}
Ef(x_1,\cdots,x_n)=
\end{equation*}
\begin{equation*}
(\frac{n}{2\pi})^{\frac{n}{2}}\int_{-\infty}^{+\infty}\cdots\int_{-\infty}^{+\infty}f(x_1,\cdots,x_n)
\exp\{-\frac{n}{2}\sum_{m=1}^n(x_m-x_{m-1})^2\}\mathrm{d}x_1\cdots\mathrm{d}x_n
\end{equation*}
\begin{equation*}
=(\frac{n}{2\pi})^{\frac{n}{2}}\int_{-\infty}^{+\infty}\cdots\int_{-\infty}^{+\infty}f(y_1,y_1+y_2,\cdots,y_1+\cdots+y_n)
\exp\{-\frac{n}{2}\sum_{m=1}^ny^2_m\}\mathrm{d}y_1\cdots\mathrm{d}y_n.
\end{equation*}
For example, taking $Y=\int_0^1x(t)\mathrm{d}t$, we have
\begin{equation*}
Y_n=\frac{1}{n}\sum_{k=1}^n x_k=\frac{1}{n}\sum_{k=1}^n(y_1+\cdots+y_k)=\frac{1}{n}(ny_1+(n-1)y_2+\cdots+y_n).
\end{equation*}
Therefore,
\begin{equation*}
EY_n=\frac{1}{n}(nEy_1+(n-1)Ey_2+\cdots+Ey_n)=0,
\end{equation*}
and
\begin{equation*}
DY_n=\frac{1}{n^2}(n^2Dy_1+(n-1)^2Dy_2+\cdots+Dy_n)
\end{equation*}
\begin{equation*}
=\frac{1}{n^3}(n^2+(n-1)^2+\cdots+1)=\frac{n(n+1)(2n+1)}{6n^3}.
\end{equation*}
So, we have
\begin{equation*}
EY=0,
\end{equation*}
and
\begin{equation*}
DY=\frac{1}{3}.
\end{equation*}

By the same method, we can compute the mean value and the variance of $Y=\int_0^1x^2(t)\mathrm{d}t$ and get the same results as the previous section.

\end{document}